\documentclass{paper1}

\usepackage{ctrfont}
\usepackage{natbib}

\usepackage{graphicx}
\usepackage{subcaption}
\usepackage{epstopdf}

\usepackage{amsmath, rotating, amsfonts}

\usepackage{url}

\usepackage{color}
\usepackage{array}

\newcolumntype{L}[1]{>{\raggedright\let\newline\\\arraybackslash\hspace{0pt}}m{#1}}
\newcolumntype{C}[1]{>{\centering\let\newline\\\arraybackslash\hspace{0pt}}m{#1}}
\newcolumntype{R}[1]{>{\raggedleft\let\newline\\\arraybackslash\hspace{0pt}}m{#1}}

\usepackage{graphicx}
\usepackage{pifont}
\usepackage{ragged2e}
\usepackage[font=sl, labelfont={sf}, margin=1cm]{caption}
\usepackage{subcaption}
\DeclareMathOperator*{\argmin}{\arg\!\min}

\usepackage[utf8]{inputenc}
\usepackage[english]{babel}
\newtheorem{theorem}{Theorem}

\usepackage{algorithm2e, algorithmic}

\title{Adjoint-based Gradient Estimation \\ Using the Space-time Solutions \\ of Unknown 
Conservation Law Simulations} 

\author{Han CHEN\footnote[1]{Aerospace Computational Design Laboratory, Massachusetts Institute of Technology} \and Qiqi WANG\footnotemark[1]}

\begin{document}
\maketitle

\begin{abstract}
Many control applications can be formulated as optimization constrained by conservation laws. 
Such optimization can be efficiently solved by gradient-based methods, where the gradient
is obtained through the adjoint method. 
Traditionally, the adjoint method has not been able to be implemented in "gray-box" conservation law simulations. 
In gray-box simulations, the analytical and numerical form of the conservation law is unknown, 
but the space-time solution of relevant flow quantities is available. 
Without the adjoint gradient, optimization can be challenging for problems with 
many control variables. However, much information about the gray-box simulation is
contained in its space-time solution, which motivates us to estimate the adjoint gradient 
by leveraging the space-time solution.
This article considers a type of gray-box simulations where the flux function is partially unknown.
A method is introduced to estimate the adjoint gradient at a cost independent of the number of control
variables. The method firstly infers a conservation law, named the twin model,
from the space-time solution, and then 
applies the adjoint method to the inferred twin model to estimate the gradient.
The method is demonstrated to achieve good gradient estimation accuracies in several numerical
examples.
The main contributions of this paper are:
a twin model method that enables the adjoint gradient computation for gray-box
conservation law simulations; and an adaptive basis construction 
scheme that fully exploits the information of gray-box solutions.
\end{abstract}

\textbf{Notations}\\
\begin{itemize}
    \item $t\in [0,T]$: the time,
    \item $\{t_i\}_{i=1}^M$: the time discretization,
    \item $x\in \Omega$: the space,
    \item $\{x_j\}_{j=1}^N$: the space discretization,
    \item $u$: the space-time solution of gray-box conservation law,
    \item $\tilde{u}$: the space-time solution of twin-model conservation law,
    \item $\boldsymbol{u}$: the discretized space-time solution of gray-box simulator,
    \item $\tilde{\boldsymbol{u}}$: the discretized space-time solution of twin-model simulator,
    \item $k$: 1) the number of equations of the conservation law; or 2) the number of folds in 
               cross validation.
    \item $D$: a differential operator,
    \item $F$: the unknown function of the gray-box model,
    \item $\tilde{F}$: the inferred $F$,
    \item $q$: the source term,
    \item $c$: the control variables,
    \item $w$: the quadrature weights in the numerical space-time integration,
    \item $\xi$: the objective function,
    \item $c_{\min}, \, c_{\max}$: bound constraints,
    \item $d$: the number of control variables,
    \item $\mathcal{C}\subset \mathbb{R}^d$: the control space,
    \item $\mathcal{M}$: the solution mismatch,
    \item $\overline{\mathcal{M}}$: the mean solution mismatch in cross validation,
    \item $\phi$: the basis functions for $\tilde{F}$,
    \item $\alpha$: the coefficients for $\phi$,
    \item $\mathcal{A}$: the index set for a basis dictionary,
    \item $T$: twin model,
    \item $\tau$: residual,
    \item $\boldsymbol{\tau}$: discretized residual,
    \item $\mathcal{T}$: integrated truncation error.
\end{itemize}

\section{Motivation}
\label{sec: motivation}
A conservation law states that a particular property of a physical system does not
appear or vanish as the system evolves over time, such as the conservation of mass, momentum, and energy. 
Mathematically, a conservation law can be expressed locally as a continuity equation
\eqref{eqn: conservation law},
\begin{equation}
    \frac{\partial u}{\partial t} + \nabla \cdot F = q\,,
    \label{eqn: conservation law}
\end{equation}
where $u$ is the conserved physical quantity, $t$ is time, $F$ is the flux of $u$, and $q$
is the source for $u$. Many equations fundamental to the physical world, such as the Navier-Stokes equation,
the Maxwell equation, and the porous medium transport equation, can be described by \eqref{eqn: conservation law}.\\

Optimization constrained by conservation laws is present in many engineering applications.
For example, in gas turbines, the rotor blades can operate at a temperature close to 2000K \cite{airfoil cooling}.
To prevent material failure due to overheating, channels can be drilled inside the rotor blades
to circulate coolant air whose dynamics are governed by the Navier-Stokes equation 
\cite{ubend rans opt 1}. The pressure
used to drive the coolant flow is provided by the compressor, resulting in a penalty on the 
turbine's thermo-dynamic efficiency \cite{ubend rans opt 2}. 
Engineers are thereby interested in optimizing the coolant
channel geometry in order to suppress the pressure loss. In this optimization problem, the control variables
are the parameters that describe the channel geometry. The dimensionality of the optimization is the
number of control variables, i.e. the control's degree of freedom.
Another example is the field control of petroleum
reservoir. In petroleum reservoir, the fluid flow of various phases and chemical components 
is dictated by porous medium transport equations \cite{reservoir sim book}. The flow can be passively and actively controlled by a variety of techniques \cite{first reservoir opt},
such as the wellbore pressure control, the polymer injection, and the steam heating,
where the reservoir is controlled by the pressure at each wells, by the the injection rate of polymer, and by the temperature of the steam \cite{secondary recovery review}.
The pressure, injection rate, and temperature can vary in each well and at every day over decades of continuous operations.
The dimensionality of the optimization is
the total number of these control variables.
Driven by economic interests, petroleum producers are 
devoted to optimizing the controls for enhanced recovery and reduced cost.\\

Such optimization is being revolutionized by the numerical simulation and optimization
algorithms. On one hand, conservation law simulation can provide an evaluation of a candidate control
that is cheaper, faster, and more scalable than conducting physical experiments. 
On the other hand, advanced optimization algorithms can guide 
the control towards the optimal with reduced number of simulation \cite{quasiNewton, gradfreereview, LBFGS, trustregionconn, genetic algo, particle swarm, cuckoo, review EI}.
However, optimization based on conservation law simulation can still be overwhelmingly costly. 
The cost is two-folded: Firstly, each simulation for a given control may run for hours or days even on a high-end
computer. This is mainly because of the high-fidelity physical models, the complex numerical schemes, and the large scale space-time 
discretization employed in the simulation. Secondly, optimization algorithms generally take many iterations of simulation on
various controls. The number of iterations required to achieve near-optimality usually 
increases with the control's degree of freedom \cite{opt via sim review}. The two costs are multiplicative.
The multiplicative effect compromises the impact of computational efforts among field engineers.\\

Fortunately, the cost due to iteration can be alleviated by adopting gradient-based optimization algorithms \cite{opt via sim review}.
A gradient-based algorithm usually 
requires significantly less iterations than a derivative-free algorithm for problems with many control
variables \cite{intro adjoint, opt via sim review, gradfreereview}.
Gradient-based algorithms require the gradient of the optimization objective to the control variables, which is efficiently 
computable through the adjoint method \cite{adjoint}. The adjoint method propogates the gradient from the objective 
backward to the control variables through the path of time integration \cite{adjoint} or through the chain of numerical operations \cite{AD review}. 
To keep track of the back propogation, the simulator source code needs to be available.
In real-world industrial simulators, adjoint is scarcely implemented because most source codes are proprietary and/or legacy.
For example, \textit{PSim}, a reservoir simulator developed and owned by \textit{ConocoPhillips}, 
is a multi-million-line Fortran-77 code that traces its birth back to the 1980's. 
Implementing adjoint directly into the source code is unpreferable because it can take tremendous amount of brain hours.
Besides, the source code and its physical models are only accessible and modifiable by the computational team inside the company. 
For the sake of gradient computation, \textit{PSim} has been superceded by adjoint-enabled simulators, but it is
difficult to be replaced due to its legacy use and cost concerns.
The proprietary and legacy nature of many industrial simulators hinders the prevalence of the adjoint method 
and gradient-based algorithms in many real-world problems with high-dimensional control.\\

Despite their proprietary and legacy nature, most simulators for unsteady conservation laws are able to 
provide the discretized space-time solution of relevant flow quantities. 
For example, \textit{PSim} provides the space-time solution of pressure, saturation, and concentration for multi-phase flow.
Similarly, most steady state simulators are able 
to provide the spatial solution. the discussion will focus on the unsteady case, 
since a steady state simulator can be viewed as a special case of the unsteady
one where the solution remains the same over many time steps.\\

I argue that the adjoint gradient computation may be enabled by leveraging the space-time solution.
The discretized space-time solution provides invaluable information about the conservation law
hardwired in the simulator. For illustration, consider a code which simulates
\begin{equation}\begin{split}
    &\frac{\partial u}{\partial t} + \frac{\partial F(u)}{\partial x} = c\,,\quad
    x \in [0,1] \;, \;\; t\in[0,1]
    \label{eqn: 1D conceptual}
\end{split}\end{equation}
with proper intial and boundary conditions and $F$ being differentiable. 
$c$ indicates the control that acts as a source for $u$.
If the expression of $F(u)$ in the simulator is not accessible by the user, 
adjoint can not be implemented directly. However, $F$ may be partially inferred from 
a discretized space-time solution of $u$ for a given $c$.
To see this, let the discretized solution be $\boldsymbol{u}\equiv \{u(t_i, x_j)\}_{i=1,\cdots,M\,,\; j=1,\cdots, N}$,
where $0\le t_1 < t_2 < \cdots < t_M \le 1$ and $0 \le x_1 < x_2 < \cdots < x_N \le 1$ indicate the time and space discretization. 
Given $\boldsymbol{u}$, the $\frac{\partial u}{\partial t}$ and $\frac{\partial u}{\partial x}$ can be sampled by finite difference.
Because \eqref{eqn: 1D conceptual} can be written as
\begin{equation}
    \frac{\partial u}{\partial t} + \frac{d F}{du} \frac{\partial u}{\partial x} = c\,,\quad
    x \in [0,1] \;, \;\; t\in[0,1]
    \label{eqn: 1D conceptual no shock}
\end{equation}
away from the shock wave, the samples of $\frac{\partial u}{\partial t}$ and $\frac{\partial u}{\partial x}$
can be plugged into \eqref{eqn: 1D conceptual no shock} to obtain samples of $\frac{dF}{du}$.
The reasoning remains intact at the shock wave, where $\frac{dF}{du}$ in \eqref{eqn: 1D conceptual no shock} is 
replaced by the finite difference form $\frac{\Delta F}{\Delta u}$ according to the Rankine-Hugoniot condition.
Based upon the sampled $\frac{dF}{du}$ and $\frac{\Delta F}{\Delta u}$, the unknown flux function $F$
can be approximated up to a constant 
for values of $u$ that appeared in the solution, 
by using indefinite integral. Let $\tilde{F}$ be the approximation for $F$.
An alternative conservation law can be proposed
\begin{equation}
    \frac{\partial \tilde{u}}{\partial t} + \frac{\partial \tilde{F}(\tilde{u})}{\partial \tilde{u}} = c\,,\quad
    x \in [0,1] \;, \;\; t\in[0,1]\,,
    \label{eqn: 1D conceptual inferred}
\end{equation}
that approximates the true but unknown conservation law \eqref{eqn: 1D conceptual}, where
$\tilde{u}$ is the solution associated with $\tilde{F}$, in the following sense: If $\tilde{F}$ and $F$ are off by a constant $a$, i.e.
$\tilde{F} = F+a$,
then $\frac{\partial F(u)}{\partial u} = \frac{\partial (F(u)+a)}{\partial u} = \frac{\partial \tilde{F}(u)}{\partial u}$;
therefore, the solutions of \eqref{eqn: 1D conceptual} 
and \eqref{eqn: 1D conceptual inferred}  to any initial value
problem will be the same.
In addition, the solutions to any adjoint equation, with any objective function, will be the same. As a result, the gradient computed 
by the adjoint equation of \eqref{eqn: 1D conceptual inferred} will be the true gradient, and therefore can drive
the optimization constrained by \eqref{eqn: 1D conceptual}.
A simulator for the approximated conservation law is named \textbf{twin model}, since it behaves as an adjoint-enabled
twin of the original simulator.
If a conservation law has a system of equations and/or has a greater-than-one spatial dimension, 
the above simple method to recover the flux function from a solution will no longer work. 
Nonetheless, much information about the flux function can be extracted from the solution. 
Given some additional information of the conservation law, one may be able to recover the unknown aspects of the flux function. 
The details of this topic are discussed in Section \ref{chapter 2}.\\

This article focuses on a class of simulators that I call \textbf{gray-box}.
A simulator is defined to be gray-box if the following two conditions are met:
\begin{enumerate}
    \item the adjoint is unavailable, and is impractical to implement into the source code.
    \item the full space-time solution of relevant flow quantities is available.
\end{enumerate}
Many industrial simulators, such as \textit{PSim}, satisfy both conditions.
In contrast, a simulator is named \textbf{open-box} if condition 1 is violated.
For example, \textit{OpenFOAM} \cite{openfoam} is an open-source fluid simulator where adjoint can be implemented directly
into its source code, so it is open-box by definition. Open-box simulators enjoy the
benifit of efficient gradient computation brought by adjoint, thereby are not within the research
scope of the article. If condition 1 is met but 2 is violated, a simulator
is named \textbf{black-box}. For example, \textit{Aspen} \cite{aspen}, an industrial chemical reactor simulator, provides neither
the adjoint nor the full space-time solution. Black-box simulators are simply calculators
for the objective function. Due to the lack of space-time solution, adjoint can not be enabled
using the twin model. 
Gray-box simulators are ubiquitous in many engineering applications. Examples are Fluent \cite{fluent} and CFX \cite{cfx} for computational fluid
dynamics, and ECLIPSE (Schlumberger), PSim (ConocoPhillips), and MORES (Shell) for petroleum reservoir simulations.
This article will only investigate gray-box simulators.\\

This article aims at estimating the adjoint gradient at a cost independent of
the number of control variables. 
Motivated by the adjoint gradient computation, a mathematical procedure for the estimation 
is developed by using the full space-time solution. In a companion paper,
I examine how the estimated gradient can faciliate a suitable 
optimization algorithm to reduce the number of iterations. \\

Instead of discussing gray-box simulators in general, this article only focuses on
simulators with partially unknown flux function, while their boundary condition,
initial condition, and the source term are known. For example, one may know that
the flux depends on certain variables, but the specific function form of such dependence is unknown.
This assumption is valid for some applications,
such as simulating a petroleum reservoir with polymer injection.
The flow in such reservoir is governed by multi-phase multi-component porous medium transport 
equations \cite{reservoir sim book}. The initial condition is usually given at the equilibrium state, the boundary
is usually described by a no-flux condition, and the source term can be modeled as controls
with given flow rate or wellbore pressure.
Usually the flux function is given by the Darcy's law.
The Darcy's law involves physical models like the permeability\footnote{The permeability quantifies the
easiness of liquids to pass through the rock.} and the viscosity\footnote{The viscosity quantifies the
internal friction of the liquid flow.}.
The mechanism through which the injected polymer modifies the rock permeability and flow viscosity can be unavailable. Thereby the flux is partially unknown.
The specific form of PDE considered in this article is given in Section \ref{sec: formulation}.
It is a future work to extend my research to more general gray-box settings where
the initial condition, boundary condition, source term, and the flux are jointly unknown.\\

\section{Problem Formulation}
\label{sec: formulation}
Consider an objective function
\begin{equation}
    \xi(\boldsymbol{u},c) = \sum_{i=1}^M \sum_{j=1}^N w_{ij} f(\boldsymbol{u}_{ij},c; t_i,x_j)
    \approx \int_0^T \int_\Omega f(u,c; t,x)d\boldsymbol{x} d t
    \label{eqn: objective prototype}
\end{equation}
where $\boldsymbol{u}$ is the discretized space-time solution of a gray-box conservation law
simulator. The spatial coordinate is $x \in \Omega$ and the time is $t\in[0,T]$.
$i=1,\cdots, M$ and $j=1,\cdots, N$ indicate the indices for the time and space discretization.
$f$ is a given function that depends on $u$, $c$, $t$, and $x$.
$w_{ij}$'s are given quadrature weights for the integration.
$c\in\mathbb{R}^d$ indicates the control variable. \\

The gray-box simulator solves the partial differential equation (PDE)
\begin{equation}\begin{split}
    \frac{\partial u}{\partial t}+ \nabla \cdot \big( D F(u) \big) = q(u,c)\,,
\end{split}
\label{eqn: govern PDE}
\end{equation}
which is a system of $k$ equation. 
The initial and boundary conditions are known.
$D$ is a known differential operator that may depend on $u$, and
$F$ is an unknown function that depends on $u$.
$q$ is a known source term that depends on $u$ and $c$.
Notice \eqref{eqn: govern PDE} degenerates to \eqref{eqn: conservation law} when $D$ equals 1.
The simulator does not have the adjoint capability,
and it is infeasible to implement the adjoint method into its source code. But the full space-time solution 
$\boldsymbol{u}$ is provided. The steady-state conservation law is a special case
of the unsteady one, so it will not be discussed separately.\\

\section{Literature Review}
\label{section: literature}

Given the background, I review the literature on derivative-free optimization and 
gradient-based optimization. In addition,
I review the adjoint method. Finally, I review methods for adaptive basis construction,
which is useful for the adaptive parameterization of a twin model.

\subsection{Review of Optimization Methods}
\label{sec: review optimization}
Optimization methods can be categorized into derivative-free and gradient-based methods 
\cite{gradfreereview}, depending on whether the gradient information is used. 
In the sequel, I review the two types of methods.\\
\subsubsection{Derivative-free Optimization}
\label{section: DFO}
Derivative-free optimization (DFO) requires only the availability of objective
function values but no gradient information \cite{gradfreereview}, thus is useful when the gradient
is unavailable, unreliable, or too expensive to obtain. 
Such methods are suitable for problems constrained by black-box simulators.\\

Depending on whether a local or global optimum is desired,
DFO methods can be categorized into local methods and global methods \cite{gradfreereview}.
Local methods seek a local optimum which is also the global optimum for convex problems. 
An important local method is the
trust-region method \cite{trust region review}. Trust-region method introduces a surrogate model that is
cheap to evaluate and presumably accurate within a trust region: an adaptive neighborhood around the current iterate \cite{trust region review}. 
At each iteration, the surrogate is optimized in a domain bounded by the trust region to generate
candidte steps for additional objective evaluations \cite{trust region review}.
The surrogates can be constructed either by interpolating the objective evaluations \cite{linear trust region, trustregionwild}, or by running a low-fidelity simulation \cite{MFO trust region,
Alexandrov trust region}.
Convergence to the objective function's optimum is guaranteed by ensuring that
the surrogate have the
same value and gradient as the objective function when the size of the trust region shrinks to zero
\cite{trustregionconn, trustregionwild}.\\

Global methods seek the global optimum. 
Example methods include the branch-and-bound search \cite{Branch and Bound}, evolution methods \cite{evolution review}, and 
Bayesian methods \cite{practical Bayesian, Locatelli, prob of improve}. 
The branch-and-bound search sequentially partitions the entire control space into a tree structure, and
determines lower and upper bounds for the optimum
\cite{Branch and Bound}. Partitions that are inferior are eliminated 
in the course of the search \cite{Branch and Bound}. The bounds are usually obtained through the assumption of the
Lipschitz continuity or statistical bounds for the objective function \cite{Branch and Bound}. 
Evolution methods maintain a population
of candidate controls, which adapts and mutates in a way that resembles natural phenomenons
such as the natural selection \cite{genetic algo, cuckoo} 
and the swarm intelligence \cite{particle swarm}.
Bayesian methods model the objective function as a random member function from a stochastic process.
At each iteration, the statistics of the stochastic process are calculated 
and the posterior, a probability measure, of the objective 
is updated using Bayesian metrics \cite{practical Bayesian, review EI}. 
The posterior is used to pick the next candidate step
that best balances the exploration of unsampled regions and the exploitation around the sampled
optimum \cite{Locatelli, jones1998, GP bandit}. \\

Because many real-world problems are non-convex, global methods are usually preferred to local methods
if the global optimum is desired \cite{gradfreereview}. Besides, DFO methods usually require
a large number of function evaluations to converge, especially when the dimension of control is large
\cite{gradfreereview}. This issue can be alleviated by incorporating the gradient information
\cite{derivative RKHS, grad coKriging, grad particle swarm, grad cuckoo}.
The details are discussed in the next subsection.\\

\subsubsection{Gradient-based Optimization}
\label{section: GBO}
Gradient-based optimization (GBO) requires the availability of the gradient values \cite{opt via sim review, nonlinear program}. 
A gradient value, if exists,
provides the optimal infinitesimal change of control variables at each iterate,
thus is useful in searching for a better control.
Similar to DFO, GBO can also be categorized into local methods and global methods \cite{opt via sim review}.
Examples of local GBO methods include the gradient descent methods 
\cite{stochastic search, backtrack line search}, the conjugate gradient 
methods \cite{linear conjugate gradient, nonlinear conjugate gradient}, and the quasi-Newton methods \cite{quasiNewton, LBFGS}. 
The gradient descent methods and
the conjugate gradient methods
choose the search step in the direction of either the gradient \cite{stochastic search, backtrack line search} or a conjugate gradient
\cite{linear conjugate gradient, nonlinear conjugate gradient}.
Quasi-Newton methods, such as the Broyden-Fletcher-Goldfarb-Shannon (BFGS) method
\cite{quasiNewton},
approximate the Hessian matrix using a series of gradient values. The approximated Hessian allows
a local quadratic approximation to the objective function which determines the search direction
and stepsize
by the Newton's method \cite{quasiNewton}. 
In addition, some local DFO methods can be enhanced to use gradient information \cite{trust region global, trust region inexact grad}. 
For instance,
in trust-region methods, the construction of local surrogates can incorporate 
gradient values if available \cite{trust region global, trust region inexact grad}.
The usage of gradient usually improves the surrogate's accuracy thus enhances the quality
of the search step, thereby reducing the required number of iterations
\cite{trust region global, trust region inexact grad}.\\

Global GBO methods search for the global optimum using gradient values \cite{opt via sim review, nonlinear program}.
Many global GBO methods can trace their development to corresponding DFO methods \cite{stogo 1, stogo 2, grad particle swarm, grad cuckoo, grad coKriging}. 
For example, the stochastic gradient-based global optimization method (StoGo) 
\cite{stogo 1, stogo 2} works by partitioning the control space and bounding the optimum in the same way 
as the branch-and-bound method \cite{Branch and Bound}. But the search in each partition is performed
by gradient-based algorithms such as BFGS \cite{quasiNewton}. 
Similarly, some gradient-based evolution methods,
such as the gradient-based particle swarm method \cite{grad particle swarm} 
and the gradient-based cuckoo search method \cite{grad cuckoo},
can be viewed as gradient variations of corresponding derivative-free counterparts \cite{particle swarm, cuckoo}.
For example, the gradient-based particle swarm method combines particle swarm algorithm with the
stochastic gradient descent method \cite{grad particle swarm}. 
The movement of each particle is dictated not only by the 
function evaluations of all particles, but also by its local gradient \cite{grad particle swarm}.\\

To achieve a desired objective value, 
GBO methods generally require much less iterations than DFO methods for problems with many 
control variables \cite{opt via sim review, nonlinear program}. GBO methods can be 
efficiently applied to 
optimization constrained by open-box simulators, because the gradient is efficiently computable by
the adjoint method \cite{adjoint, opt via sim review}.
In a companion paper, I extend GBO to optimization constrained by gray-box simulation by estimating 
the gradient using the full space-time solution.

\subsection{The Adjoint Method}
\label{sec: adjoint}
Consider a differentiable objective function constrained by a conservation law PDE 
\eqref{eqn: govern PDE}.
Let the objective function be $\xi(u,c)$, $c\in \mathbb{R}^d$, 
and let the PDE \eqref{eqn: govern PDE} be abstracted
as $\mathcal{F}(u,c) = 0$. 
$\mathcal{F}$ is a parameterized differential operator, 
together with boundary conditions and/or initial conditions, that uniquely defines a $u$ for each $c$.
The gradient $\frac{d\xi}{dc}$ can be estimated trivially by finite difference. The
$i$th component of the gradient is given by
\begin{equation}
    \left(\frac{d\xi}{dc}\right)_i \approx \frac{1}{\delta} \big( 
    \xi(u+ \Delta u_i, c+\delta e_i) - \xi(u, c) \big)\,,
    \label{eqn: FD eg}
\end{equation}
where
\begin{equation}\begin{split}
    \mathcal{F}(u,c)=0\,, \quad \mathcal{F}(u+\Delta u_i, c+ \delta e_i) = 0\,.
\end{split}
\label{eqn: define perturb}
\end{equation}
$e_i$ indicates the $i$th unit Cartesian basis vector in $\mathbb{R}^d$, and $\delta>0$ indicates
a small perturbation. 
Because \eqref{eqn: define perturb} needs to be solved for every $\delta e_i$, 
so that the corresponding $\Delta u_i$ can be used in \eqref{eqn: FD eg},
$d+1$ PDE simulations are required to evaluate the gradient.
As explained in Section \ref{sec: review optimization}, $d$ can be large in many
control optimization problems. Therefore, it can be costly to evaluate the gradient 
by finite difference.\\

In contrast, the adjoint method evaluates the gradient using only one PDE simulation
plus one adjoint simulation
\cite{adjoint}.
To see this, linearize $\mathcal{F}(u,c)=0$ into a variational form
\begin{equation}
    \delta \mathcal{F} = \frac{\partial \mathcal{F}}{\partial u} \delta u + \frac{\partial \mathcal{F}}{\partial c} \delta c = 0\,,
\label{eqn: variational form}
\end{equation}
which gives
\begin{equation}
    \frac{du}{dc} = - \left(\frac{\partial \mathcal{F}}{\partial u}\right)^{-1} \frac{\partial \mathcal{F}}{\partial c}
    \label{eqn: variational PDE}
\end{equation}
Using \eqref{eqn: variational PDE}, $\frac{d\xi}{dc}$ can be expressed by
\begin{equation}\begin{split}
    \frac{d\xi}{dc} &= \frac{\partial \xi}{\partial u}\frac{du}{dc} + \frac{\partial \xi}{\partial c}\\
    &= - \frac{\partial \xi}{\partial u} \left(\frac{\partial \mathcal{F}}{\partial u}\right)^{-1}
       \frac{\partial \mathcal{F}}{\partial c} + \frac{\partial \xi}{\partial c}\\
    &= - \lambda^T \frac{\partial \mathcal{F}}{\partial c} + \frac{\partial \xi}{\partial c}
\end{split}
\,,
    \label{eqn: dxidc adjoint}
\end{equation}
where $\lambda$, the adjoint state, is given by the adjoint equation
\begin{equation}
    \left(\frac{\partial \mathcal{F}}{\partial u}\right)^T \lambda = \left(\frac{\partial \xi}{\partial u}\right)^T
    \label{eqn: adjoint}
\end{equation}
Therefore, the gradient can be evaluated by \eqref{eqn: dxidc adjoint}
using one simulation of $\mathcal{F}(u,c)=0$ and one simulation of \eqref{eqn: adjoint} that 
solves for $\lambda$. \\

Adjoint methods can be categorized into continuous adjoint and discrete adjoint methods,
depending on whether the linearization or the discretization is excuted first \cite{review adjoint geo}.
The above procedure, \eqref{eqn: variational form} thru. \eqref{eqn: adjoint},
is the continuous adjoint,
where $\mathcal{F}$ is a differential operator.
The continous adjoint method linearizes the continuous PDE $\mathcal{F}(u,c) = 0$ first,
then discretizes the adjoint equation 
\eqref{eqn: adjoint} \cite{adjoint}. 
In \eqref{eqn: adjoint}, 
$\left(\frac{\partial \mathcal{F}}{\partial u}\right)^T$ can
be derived as another differential operator. With proper boundary and/or initial
conditions, it uniquely determines the adjoint solution $\lambda$. See \cite{intro adjoint}
for a detailed derivation of the continuous adjoint equation.\\

The discrete adjoint method \cite{discrete adjoint}
discretizes $\mathcal{F}(u,c)=0$ first. After the discretization, $u$ and $c$ become
vectors $\boldsymbol{u}$ and $\boldsymbol{c}$.
$\boldsymbol{u}$ is defined implicitly by the system
$\mathcal{F}_d(\boldsymbol{u}, \boldsymbol{c}) = 0$, where $\mathcal{F}_d$ indicates
the discretized difference operator. Using the same derivation as \eqref{eqn: variational form}
thru. \eqref{eqn: adjoint}, the discrete adjoint equation can be obtained
\begin{equation}
    \left(\frac{\partial \mathcal{F}_d}{\partial \boldsymbol{u}}\right)^T \boldsymbol{\lambda} = \left(\frac{\partial \xi}{\partial \boldsymbol{u}}\right)^T\,,
    \label{eqn: adjoint discrete}
\end{equation}
which is a linear system of equations.
$\left(\frac{\partial \mathcal{F}_d}{\partial \boldsymbol{u}}\right)^T$ is derived as
another difference operator. With proper discretized boundary/initial conditions, it uniquely
determines the discrete adjoint vector $\boldsymbol{\lambda}$, which subsequently determines the
gradient
\begin{equation}
    \frac{d\xi}{d\boldsymbol{c}} = - \boldsymbol{\lambda}^T \frac{\partial \mathcal{F}_d}{\partial \boldsymbol{c}} + \frac{\partial \xi}{\partial \boldsymbol{c}}
\,.
    \label{eqn: dxidc adjoint discrete}
\end{equation}
See Chapter 1 of \cite{discrete adjoint phd} for a detailed derivation of the discrete adjoint.\\

\begin{figure}
    \begin{center}
    \includegraphics[width=4cm]{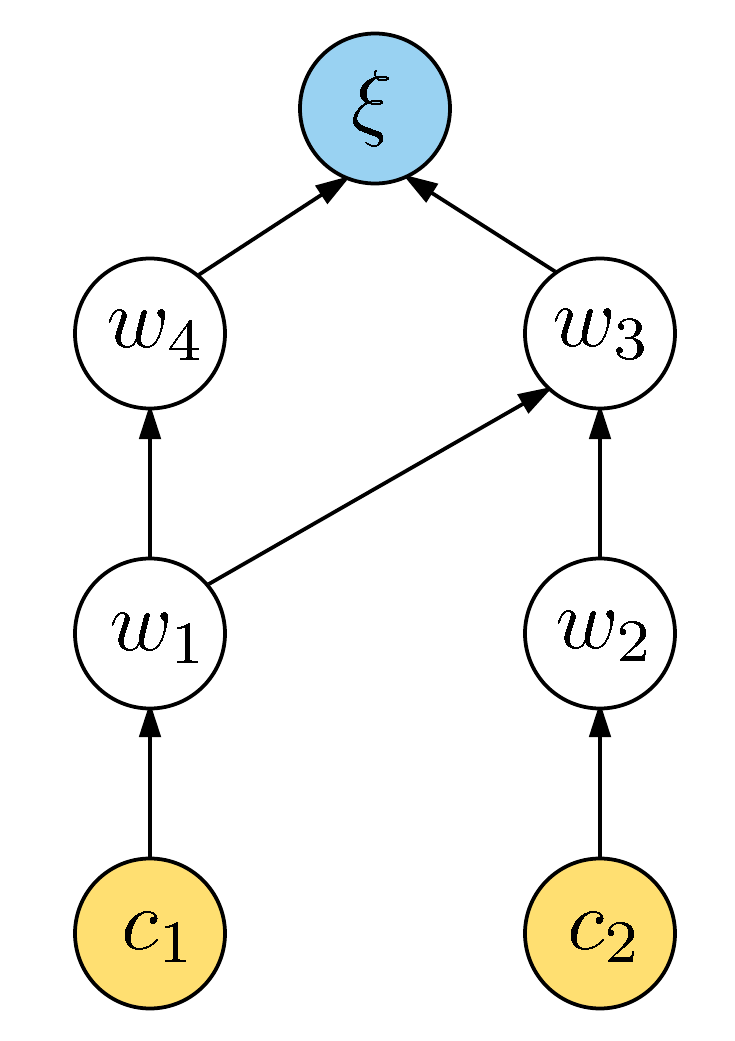}
    \caption{The computational graph for \eqref{eqn: elementary operation}. The yellow nodes indicate
             the input variables, the blue node indicates the output variable, and the white
             nodes indicate the intermediate variables. The arrows indicate elementary operations.
             The begining and end nodes of each 
             arrow indicate the independent and dependent variables for each operation.}
    \label{fig: AD graph}
    \end{center}
\end{figure}

The discrete adjoint method can be implemented by automatic differentiation (AD) 
\cite{AD review}.
AD exploits the fact that a PDE simulation, no matter how complicated, executes a sequence
of elementary arithmetic operations (e.g. addition, multiplication) and elementary functions
(e.g. $\exp$, $\sin$) \cite{AD review}. For example, consider the function
\begin{equation}
    \xi = f(c_1, c_2) = c_1 c_2 + \sin(c_1)\,.
    \label{eqn: eg y=fx}
\end{equation}
The function can be broken down into a series of elementary arithmetic operations and
elementary functions.
\begin{equation}\begin{split}
    w_1 &= c_1\\
    w_2 &= c_2\\
    w_3 &= w_1 w_2\\
    w_4 &= \sin(w_1)\\
    \xi &= w_3 + w_4\;.
    \end{split}
    \label{eqn: elementary operation}
\end{equation}
\eqref{eqn: elementary operation} can be represented by a computational 
graph in Figure \ref{fig: AD graph}.
In the graph,
the gradient of the output with respect to the input variables 
can be computed using the chain rule \cite{AD review}. 
Let $\bar{z}$ 
denote the gradient of $\xi$ with respect to $z$, for any independent or intermediate
variable $z$ in \eqref{eqn: elementary operation}.
To compute the derivatives $\bar{c}_1 = \frac{\partial \xi}{\partial c_1}$ and
$\bar{c}_2 = \frac{\partial \xi}{\partial c_2}$,
one can propogate the derivatives backward in the computational graph as follows
\begin{equation}\begin{split}
    \bar{w}_4 & = 1\\
    \bar{w}_3& =1\\
    \bar{w}_2 & = \bar{w}_3 \frac{\partial w_3}{\partial w_2} = 1\cdot w_1\\
    \bar{w}_1 & = \bar{w}_4 \frac{\partial w_4}{\partial w_1} + \bar{w_3} 
                  \frac{\partial w_3}{\partial w_1} = 1\cdot \cos(w_1) +1\cdot w_2\\
    \bar{c}_2 &= \bar{w}_2 = c_1\\
    \bar{c}_1 &= \bar{w}_1 = \cos(c_1) + c_2
\end{split}
\label{eqn: elementary derivative}
\end{equation}
The derivatives in \eqref{eqn: elementary derivative} are straightforward to compute.
This is because every forward operation in \eqref{eqn: elementary operation} is elementary, and
their derivatives can be hardwired in AD softwares. 
Notice each arrow in Figure \ref{fig: AD graph} is traversed once and only once in the backward propogation 
\eqref{eqn: elementary derivative}. Therefore, the backward
gradient computation has a similar cost as the forward output computation, regardless of
the number of input variables. 
See \cite{AD review} for a thorough review of AD.\\

\begin{figure}\begin{center}
    \begin{subfigure}[t]{.99\textwidth}
        \centering
        \includegraphics[width=7cm]{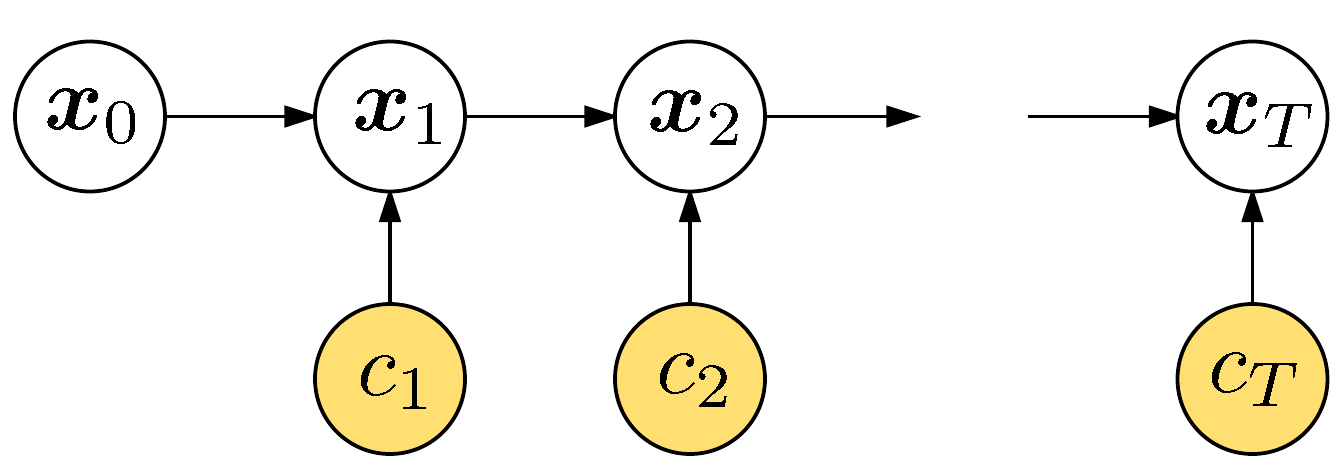}
        \caption{The computational graph for \eqref{eqn: residual PDE eg}, which is constructed
                 by \eqref{eqn: elements for markov graph}. The yellow
                 nodes indicate the input variables.}
        \label{fig: chap1 graph 1}
    \end{subfigure}
    \begin{subfigure}[t]{.99\textwidth}
        \centering
        \includegraphics[width=7cm]{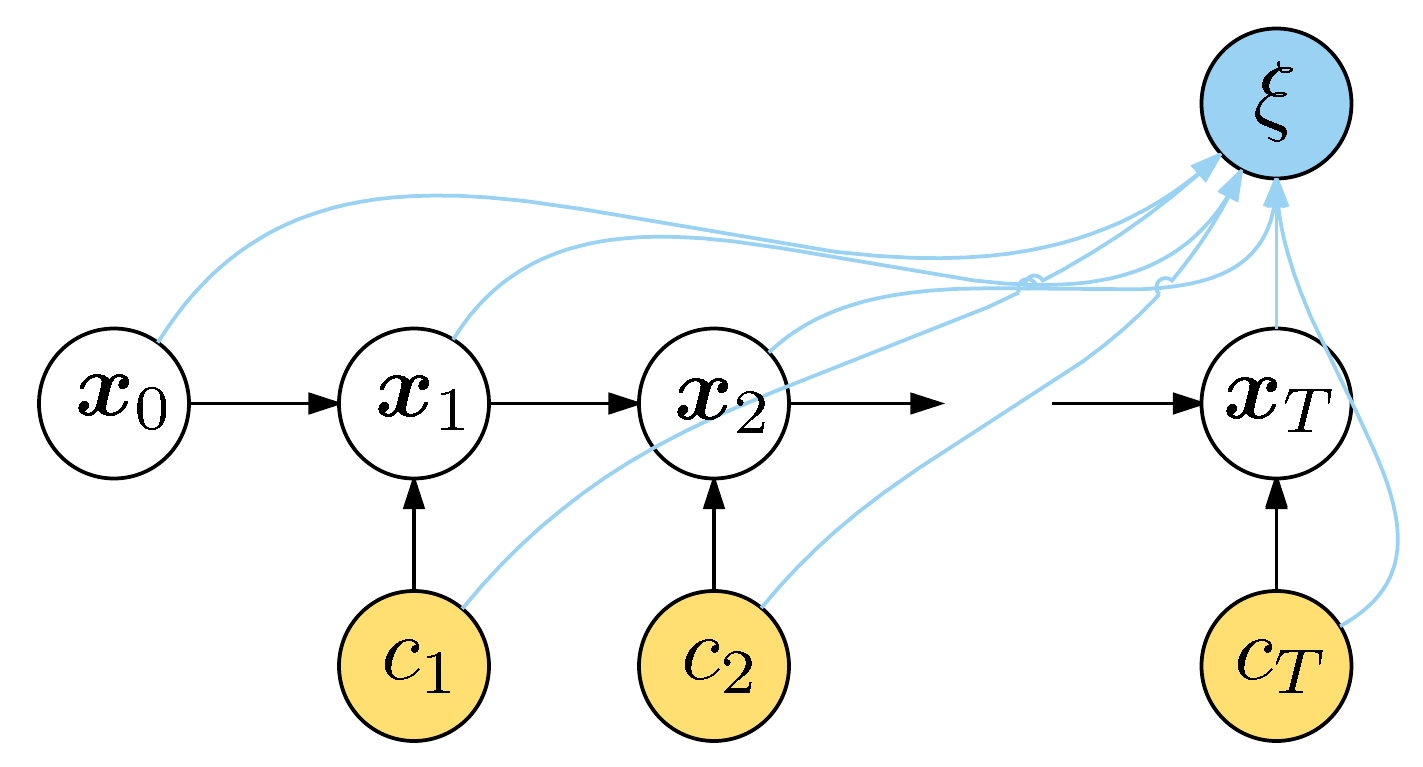}
        \caption{The computational graph for evaluating the objective function $\xi$.
                 The blue node indicates the output variable.}
        \label{fig: chap1 graph 2}
    \end{subfigure}
    \caption{Computational graphs for the PDE simulation and objective evaluation.}
\end{center}\end{figure}

Because a PDE simulation can be viewed as 
performing a sequence of elementary operations, 
AD can be used to evaluate the discrete adjoint.
Consider a discretized PDE simulation that at each timestep solves
\begin{equation}
    \mathcal{F}_{t+1} = \mathcal{F}(\boldsymbol{x}_t, \boldsymbol{x}_{t+1}, \boldsymbol{c}_{t+1})=0\,,
    \label{eqn: residual PDE eg}
\end{equation}
for $t=0, \cdots, T-1$, where $\boldsymbol{x}_t$ and $\boldsymbol{c}_t$ are
the state and control variables at the $t$th timestep. 
AD can be used to compute the gradient of an objective function
$$\xi = \xi(\boldsymbol{x}_0 , \cdots, \boldsymbol{x}_T; \boldsymbol{c}_1, \cdots 
\boldsymbol{c}_{T})$$ 
to the control variables.
To see this, consider the evaluation of \eqref{eqn: residual PDE eg} using an AD software. 
The gradients $\frac{\partial \mathcal{F}_{t+1}}{\partial \boldsymbol{x}_t}$, 
$\frac{\partial \mathcal{F}_{t+1}}{\partial \boldsymbol{x}_{t+1}}$, and
$\frac{\partial \mathcal{F}_{t+1}}{\partial \boldsymbol{c}_{t+1}}$, for
$t=0, \cdots, T-1$, are automatically computable. Therefore, one can obtain
\begin{equation}
    \begin{split}
        \frac{\partial \boldsymbol{x}_{t+1}}{\partial \boldsymbol{x}_t} &=
        - \left(\frac{\partial \mathcal{F}_{t+1}}{\partial \boldsymbol{x}_{t+1}}\right)^{-1}
        \left(\frac{\partial \mathcal{F}_{t+1}}{\partial \boldsymbol{x}_t}\right)\\
        \frac{\partial \boldsymbol{x}_{t+1}}{\partial \boldsymbol{c}_{t+1}} &=
        -\left(\frac{\partial \mathcal{F}_{t+1}}{\partial \boldsymbol{x}_{t+1}}\right)^{-1}
        \left(\frac{\partial \mathcal{F}_{t+1}}{\partial\boldsymbol{c}_{t+1}}\right)\,.
    \end{split}
    \label{eqn: elements for markov graph}
\end{equation}
Therefore a computational graph, Figure \ref{fig: chap1 graph 1},
can be constructed using the chain rule. The graph enables the evaluation of all
$\frac{\partial \boldsymbol{x}_t}{\partial \boldsymbol{c}_{t-i}}$, for $t=1, \cdots, T$
and $i=0, \cdots, t-1$.
Given the solutions $\boldsymbol{x}_t$'s and the controls $\boldsymbol{c}_t$'s, 
the evaluation of $\xi$ is nothing but overlaying the graph by an additional 
layer of computations, shown in Figure \ref{fig: chap1 graph 2}.
Because $\frac{\partial \xi}{\partial \boldsymbol{x}_t}$'s and 
$\frac{\partial \xi}{\partial \boldsymbol{c}_t}$'s can be obtained by AD, the gradient
\begin{equation} 
    \frac{d\xi}{d \boldsymbol{c}_t}= \frac{\partial \xi}{\partial \boldsymbol{c}_t} + 
    \frac{\partial \xi}{\partial \boldsymbol{x}_t}\frac{\partial \boldsymbol{x}_t}{\partial 
       \boldsymbol{c}_t} +
    \frac{\partial \xi}{\partial \boldsymbol{x}_{t+1}} \frac{\partial \boldsymbol{x}_{t+1}}{\partial
       \boldsymbol{c}_t} + \cdots +
    \frac{\partial \xi}{\partial \boldsymbol{x}_{T}} \frac{\partial \boldsymbol{x}_{T}}{\partial
       \boldsymbol{c}_t}
\end{equation}
can be computed, for all $t=1, \cdots , T$.\\

The adjoint method has seen wide applications in optimization problems constrained by 
conservation law simulations, such as in
airfoil design \cite{adjoint aerodynamics, adjoint aerodynamics 2, 
adjoint aerodynamics AD}, adaptive mesh refinement
\cite{discrete adjoint phd}, injection policy optimization in petroleum reservoirs
\cite{adjoint reservoir optimal control}, 
history matching in reservoir geophysics
\cite{review adjoint geo}, and optimal well placement in reservoir management 
\cite{adjoint well place}.
Besides, there are many free AD softwares available for various languages, such as
\emph{ADOL-C} (C, C++) \cite{adolc}, \emph{Adiff} (Matlab) \cite{adiff}, and \emph{Theano} (Python)
\cite{theano}.
Unfortunately, the adjoint method is not directly applicable
to gray-box simulations, as explained in Section \ref{sec: motivation}. 
To break this limitation, section \ref{chapter 2} develops the twin model method
that enables the adjoint gradient computation for gray-box simulations.\\

\subsection{Adaptive Basis Construction}
\label{sec: adaptive basis review}

The unknown function $F$ in \eqref{eqn: govern PDE} can be approximated by a linear
combination of basis functions. An over-complete or
incomplete set of bases can negatively affect the approximation due to overfitting or underfitting.
Therefore, adaptive basis construction is needed.\\

Square-integrable functions
can be represented by the parameterization\footnote{
    By definition, the square 
    integrable functions form a pre-Hilbert space with inner product given by
    $<f,g> = \int_A \overline{f(x)}g(x) dx\,,$ where 
    1) $f$ and $g$ are square integrable functions, 2) $\overline{f(x)}$ is 
    the complex conjugate of $f(x)$, and 3) $A$ is the set over which the integral is defined.
    It can be shown that square integrable functions are complete 
    under the metric induced by the inner product, thus is a Hilbert space \cite{functional analysis}.
    Any member function can be represented by 
    a finite or countably infinite number of basis functions of the Hilbert space
    \cite{functional analysis}.
} \cite{functional analysis}
\begin{equation}
    \tilde{F}(\cdot) = \sum_{i\in \mathbb{N}} {\alpha_i} \phi_i(\cdot)\,,
    \label{eqn: linear param general}
\end{equation}
where $\phi_i$'s are linearly-independent basis functions, $\alpha_i$'s are the coefficients, and $i$
indices the basis. For example, 
a bivariate function can be represented by monomials (Weierstrass approximation 
theorem \cite{functional analysis})
\begin{equation*}
    1,\,  u_1,\, u_1^2,\,
    u_2,\, u_1 u_2,\,  u_1^2 u_2, \,
    u_2^2,\,  u_1 u_2^2,\,  u_1^2 u_2^2,\, \cdots\,.
\end{equation*}
on any real interval $[a,b]$.\\

Let $\mathcal{A}$ be a non-empty finite subset of $\mathbb{N}$,
$\tilde{F}$ can be approximated using a subset of bases, 
\begin{equation}
    \tilde{F}(\cdot) \approx \sum_{i\in \mathcal{A}} {\alpha_i} \phi_i(\cdot)\,,
    \label{eqn: linear param truncate}
\end{equation}
where $\{\phi_i\}_{i\in \mathcal{A}}$ is called a basis dictionary \cite{match pursuit}. 
The approximation is solely determined by the choices of the dictionary and the coefficients.
For example, in polynomial approximation, the basis dictionary can consist of
the basis whose total polynomial degree does not exceed $p\in \mathbb{N}$ \cite{PCE}.
Given a dictionary, the coefficients for $\tilde{F}$
can be determined by the minimization \cite{PCE}
\begin{equation}\begin{split}
    \boldsymbol{\alpha}^* = \argmin_{\boldsymbol{\alpha}\in \mathbb{R}^{|\mathcal{A}|}}
    \left\| \tilde{F}-\sum_{i\in \mathcal{A}} \alpha_i \phi_i \right\|_{L_p}
\end{split}\,,
\label{eqn: least squares demo}
\end{equation}
where $\|\cdot\|_{L_p}$ indicates the $L_p$ norm\footnote{
    Usually $p=1$ \cite{L1 basis pursuit} or $2$ \cite{L2 frame, match pursuit}.
}.
This article parameterizes the twin-model flux $\tilde{F}$ and
optimizes the coefficients,
so the twin model serves as a proxy of the gray-box model. Details are discussed in 
Section \ref{sec: flux param}.\\

If the dictionary is pre-determined without using any evaluation of the underlying function,
its cardinality can increase as the number of variables increases, and
as the basis complexity increases \cite{PCE}.
For example, for $d$-variate polynomial basis, the total number of bases 
is $d^p$ if one bounds the polynomial degree of each variable by $p$; and 
is $\binom{p+d}{d}$ if one bounds the total degree by $p$ \cite{PCE}.
For piecewise linear basis,
one can approximate the function using
Smolyak's sparse grid with
$\mathcal{O}\left( n(\log n)^{d-1} \right)$ bases\footnote{
To be precise, if the underlying function has bounded mixed second derivatives
$D^{\boldsymbol{\beta}} \cdot = \frac{\partial^{|\boldsymbol{\beta}|_1}\cdot }{\partial
u_1^{\beta_1} \cdots \partial u_d^{\beta_d}}$ for $|\boldsymbol{\beta}|_\infty\le 2$,
and is defined in a $d$-dimensional unit cube, then the error between
the underlying function $\tilde{F}$ and the approximation $\hat{\tilde{F}}$ is 
$\left\| \tilde{F} - \hat{\tilde{F}} \right\|_{L_2} = \mathcal{O}(4^{-n} n^{d-1})$
\cite{Smolyak}.}, where $n$ is the number of univariate basis for each variable.\\

In many applications, 
one may deliver a similarly accurate approximation by 
using a much smaller subset of the dictionary as the bases
than using the full dictionary \cite{PCE, L1 basis pursuit, match pursuit, Lasso variable selection}.
To exploit the sparse structure, only significant bases shall be selected, and
the selection process shall be adaptive depending on the values of function evaluations.
There are several methods that adaptively determine the sparsity, such as Lasso regularization
\cite{Lasso variable selection}, matching 
pursuit \cite{match pursuit}, and basis pursuit \cite{L1 basis pursuit}. 
Lasso regularization adds a penalty $\lambda
\sum_{i\in\mathcal{A}}|\alpha_i|$
to the approximation error, where $\lambda>0$ is a tunable parameter \cite{Lasso variable selection}. 
In this way, Lasso balances the approximation error and the number of non-zero
coefficients \cite{Lasso variable selection}.
Matching pursuit adopts a greedy, stepwise approach \cite{match pursuit}.
It either selects a significant basis one-at-a-time (forward selection) from a dictionary
\cite{forward selection},
or prunes an insignificant basis one-at-a-time (backward pruning) from the dictionary
\cite{backward prune}. 
Basis pursuit minimizes $\|\boldsymbol{\alpha}\|_{L_1}$ subject to \eqref{eqn: linear param general},
which is equivalently reformulated and efficiently solved as a linear programming problem 
\cite{L1 basis pursuit}.\\

Conventionally, the dictionary for the sparse approximation needs to
be pre-determined, with the belief that the dictionary is a superset of the significant bases
\cite{adaptive basis 2}. 
This can be problematic because 
the maximum complexity\footnote{
The definition of complexity is basis-dependent. For example, the complexity for
polynomial basis can be its total polynomial degree; and the complexity for
wavelet basis can be its finest resolution \cite{wavelet mallat}. Here the ``complexity'' is
discussed in a general sense. } 
of the significant ones are unknown a prior.
To address this issue, methods have been devised that construct an adaptive dictionary
\cite{adaptive basis 1, adaptive basis 2, adaptive basis 3}.
Although different in details, such methods share the same approach: Starting from some trivial
bases\footnote{
For example, the starting basis can be $1$ for polynomial basis \cite{adaptive basis 1}.
The starting bases serve as seeds from which more complex bases grow. See \cite{adaptive basis 1,
adaptive basis 2, adaptive basis 3} for more details of the heuristics. The problem will be 
revisited in Section \ref{sec: adaptive basis} and \ref{sec: twin algo}.}, 
a dictionary is built up progressively 
by iterating over a forward step and a backward step
\cite{adaptive basis 1, adaptive basis 2, adaptive basis 3}.
The forward step searches over a candidate set of bases, and appends the significant ones 
to the dictionary
\cite{adaptive basis 1, adaptive basis 2, adaptive basis 3}.
The backward step searches over the current dictionary, and removes the 
insignificant ones from the dictionary
\cite{adaptive basis 1, adaptive basis 2, adaptive basis 3}.
The iteration stops only when no alternation is made to the dictionary or when a targeted 
accuracy is achieved, without bounding the basis complexity a prior
\cite{adaptive basis 1, adaptive basis 2, adaptive basis 3}.
Such approach is adopted to build up 
the bases for $\tilde{F}$. Details are discussed in Section \ref{sec: adaptive basis}.\\

Based the motiviation and literature review, we find a need to enable adjoint
gradient computation for gray-box conservation law simulations,
especially for problems with many control variables. 
The objective is
to develop an adjoint approach that estimates the gradient of objective functions
constrained by gray-box conservation law simulations with partially unknown flux functions,
by leveraging the space-time solution.
Section \ref{infer} 
devises a general framework to estimate the gradient of an
objective function constrained by a gray-box simulation, 
at a cost independent of the gradient's dimensionality. 
This is achieved through firstly training a twin model, 
then applying the adjoint method to the trained twin model.
Section \ref{sec: flux param} introduces a parameterization of the unknown flux function,
followed by a numerical demonstration that motivates the needs for an adaptive parameterization.
Section \ref{sec: adaptive basis}
develops an adaptive scheme for the flux function.
By summarizing the developments, a twin model algorithm is presented
in section \ref{sec: twin algo}. To reduce the computational cost in training a twin model,
a truncation error metric and a pre-train step is developed in section \ref{sec: trunc error}.
Finally, Section \ref{sec: twin numerical results} demonstrates the twin model algorithm 
in several numerical examples.

\section{Estimate the Gradient by Using the Space-time Solution}
\label{chapter 2}
This section develops a method to estimate the gradient by using the space-time solution
of gray-box conservation law simulations.
An example, equation \eqref{eqn: conservation law}, has
been given in Section \ref{sec: motivation} to illustrate why such estimation is possible
when the conservation law involves only one equation and one dimensional space.
In this example, the derivative of the flux can be first interpolated by using the
discretized gray-box solution;
then the adjoint method is applied to the inferred conservation law \eqref{eqn: 1D conceptual inferred}
to estimate the gradient. This section develops a more general
procedure suitable for systems of equations and for problems with a spatial dimension greater than one.

\subsection{Approach}
\label{infer}
Consider a gray-box simulator that solves the PDE \eqref{eqn: govern PDE}, 
a system of $k$ equations, for $u(t,x)$ with $t\in[0,T]$
and $x\in\Omega$. The PDE has an unknown flux $F$, but known source term $q$, and
known initial and boundary conditions.
Let its discretized space-time solution be $\boldsymbol{u}$.
The article introduces an open-box simulator solving another PDE, namely the twin model,
\begin{equation}\begin{split}
    \frac{\partial \tilde{u}}{\partial t}+ \nabla \cdot \big(D \tilde{F}(\tilde{u})\big) = q(\tilde{u},c)\,,
\end{split}
\label{eqn: govern twin model}
\end{equation}
which is a system of $k$ equations
with the same source term and the same initial and boundary conditions. 
Equation \eqref{eqn: govern twin model}
differs from \eqref{eqn: govern PDE} in its flux. For simplicity, 
let the open-box simulator use the same space-time
discretization, and let its discretized solution be $\tilde{\boldsymbol{u}}$.
Define the solution mismatch
\begin{equation}
    \mathcal{M}(\tilde{F}) = \sum_{i=1}^M \sum_{j=1}^N w_{ij} \big( \tilde{\boldsymbol{u}}_{ij}
     -\boldsymbol{u}_{ij}\big)^2 \,,
    \label{eqn: solution mismatch}
\end{equation}
where $w_{ij}$'s are the quadrature weights for the space-time integration.
$\mathcal{M}$ approximates the space-time integration of the continous solutions' mismatch,\\
$\mathcal{M} \approx \int_{0}^T\int_\Omega \big(\tilde{u}(t,x) - u(t,x)\big)^2 \, \textrm{d}x\,\textrm{d}t$.
Given a set $\mathcal{S}_F$ consisted of all possible guesses for $F$,
I propose to infer a flux $\tilde{F}$ such that the mismatch between
$\boldsymbol{u}$ and $\tilde{\boldsymbol{u}}$ is minimized, i.e.
\begin{equation}
    \tilde{F}^* = \argmin_{\tilde{F} \in \mathcal{S}_F} \mathcal{M}\,,
    \label{eqn: twin min mismatch}
\end{equation}
The choice for $\mathcal{S}_F$ will be discussed later in Section \ref{sec: flux param} and
\ref{sec: adaptive basis}. Because the twin model is open-box, \eqref{eqn: twin min mismatch}
can be solved by gradient-based methods. Once $\mathcal{M}$ is minimized, the adjoint
method can be applied to the twin model to estimate the gradient.\\

\begin{figure}\begin{center}
    \includegraphics[height=3.4cm]{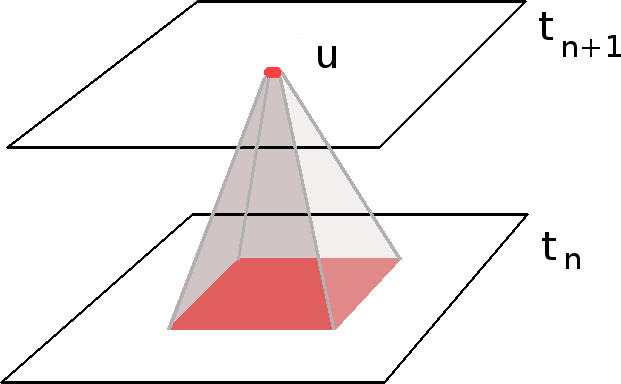}
    \caption{Domain of dependence: 
             $u(t,x)$
             depends on $u$ at an earlier time within a
             domain of dependence. The two planes in this figure indicates the spatial 
             solution at two adjacent timesteps.
             The domain of dependence can be much smaller
             than $\Omega$, the entire spatial domain.}
    \label{fig: locality}
\end{center}\end{figure}

The key to inferring the flux is to leverage the gray-box space-time solution.
Its inferrability can be 
loosely explained by the following reasonings. Firstly,
the conserved quantity $u$ in \eqref{eqn: govern PDE} depends on
$u$ in a previous time only inside a domain of dependence, illustrated in Figure \ref{fig: locality}.
Similarly, in discretized PDE simulation, the solution at any gridpoint only depends on a
numerical domain of dependence.
Besides, in a PDE simulation, a one-step time marching at any gridpoint can be viewed as
a mapping whose input only involves the numerical domain of dependence.
Because the number of state variables in the numerical domain of dependence can be small, 
inference of the mapping is potentially feasible.
Secondly, the space-time solution at every space-time gridpoint can be viewed
as a sample for this mapping. 
Because the scale of space-time discretization in a conservation law simulation is usually large,
a large number of samples are available for the inference, thus
making the inference potentially accurate.\\

\begin{figure}
    \begin{center}
        \includegraphics[width=7cm]{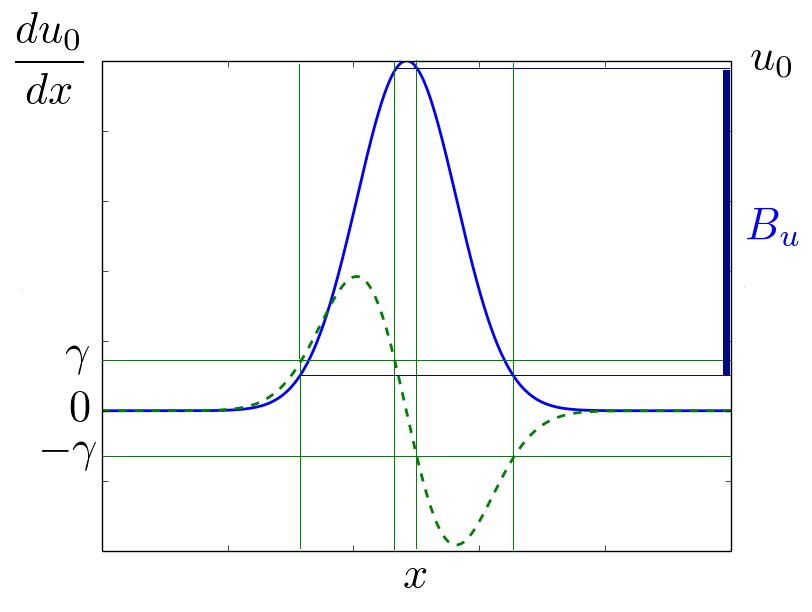}
        \caption{An illustration of $B_u$ defined in Theorem \ref{theorem: 1}. The blue line is $u_0$ 
        and the green dashed line is $\frac{du_0}{dx}$. $B_u$ is the set of $u_0$ where the derivative 
        $\frac{du_0}{dx}$ has an absolute value larger than $\gamma$.}
        \label{fig: excitedDomain}
    \end{center}
\end{figure}

The inferrability is difficult to prove.
However, it can be partially justified by the following theorem
for PDEs with $k=1$ and one dimensional space.
\begin{theorem}
    Consider two PDEs
    \begin{equation}
        \qquad\frac{\partial u}{\partial t} + \frac{\partial F(u)}{\partial x} = 0\,,\; \emph{and}
        \label{eqn: easy 1}
    \end{equation}
    \begin{equation}
        \frac{\partial \tilde{u}}{\partial t} + \frac{\partial \tilde{F}(\tilde{u})}{\partial x} = 0\,,
        \label{eqn: easy 2}
    \end{equation}
    with the same initial condition $u(0,x) = u_0(x)$, and $x\in \mathbb{R}$. $u_0$ is bounded, differentiable, 
    Lipschitz continuous with constant $L_u$, 
    and has a finite support. $F$ and $\tilde{F}$ are both twice-differentiable and Lipschtiz 
    continuous with constant $L_F$.
    Let 
    $$
    B_u \equiv \left\{ u\left| u=u_0(x) \; \emph{that satisfies} 
    \;\left|\frac{du_0}{dx}\right|\ge \gamma > 0
    \,,\; \emph{for all} \; x\in \mathbb{R}\right.
    \right\}
    \subseteq \mathbb{R}\,.
    $$
    be a non-empty and measurable set.
    We have:\\
    For any $\epsilon >0$, there exist $\delta>0$ and $T>0$ such that 
    \begin{itemize}
        \item if $|\tilde{u}(t,x)-u(t,x)| < \delta$ for any $x\in \mathbb{R}$ and $t\in [0,T]$, then
              $\left|\frac{d\tilde{F}}{du} - \frac{dF}{du}\right| < 
               \epsilon $ for any $u\in B_u\,.$
    \end{itemize}
    \label{theorem: 1}
\end{theorem}
The proof is given in Appendix \ref{proof 1}.
An illustration of $B_u$ is given in Figure \ref{fig: excitedDomain}.
Several observations can be made from Theorem \ref{theorem: 1}.
Firstly, if the solutions of \eqref{eqn: easy 1} and \eqref{eqn: easy 2} match closely, then
the derivatives of their flux functions must match closely in $B_u$. 
Secondly, the conclusion 
can only be drawn for $u\in B_u$ where the initial condition has coverage and has large enough
variation. For more general problems involving 1) systems of equations, 2) higher spatial dimensions, 
and 3) discretization,
the inferrability is difficult to show theoretically. Instead, it will be demonstrated numerically.\\

The next section discusses the choices of $\mathcal{S}_F$ occurred
in \eqref{eqn: twin min mismatch}.
A suitable parameterization for $\tilde{F}$ will be chosen that takes into account the observations from
Theorem \ref{theorem: 1}.\\

\subsection{Parameterization}
\label{sec: flux param}
Functions can be parameterized by a linear combination of basis functions \cite{functional analysis}.
Firstly, consider the case when $\tilde{F}$ is univariate. 
There are many types of basis
functions to parameterize a univariate function, such as polynomial basis, Fourier basis, and
wavelet basis. 
Based on the observations from Theorem \ref{theorem: 1}, $\tilde{F}$ and $F$
are expected to match only on a domain of $u$ 
where the gray-box space-time solution exists and has large variation. Besides, $\tilde{F}$ may match 
$F$ better on a domain where the gray-box discretized solution $\boldsymbol{u}$ are more densely sampled.
Therefore, an ideal parameterization should 
admit local refinements so $\tilde{F}$ can match $F$ better at some domain;
similarly, it should allow local dropouts when bases are redundant at some domain. 
This section presents a choice of the parameterization for $\tilde{F}$ that allows such 
local refinements, while a procedure for basis refinement is given in the next section.\\

\begin{figure}\begin{center}
    \includegraphics[width=4cm, height=3cm]{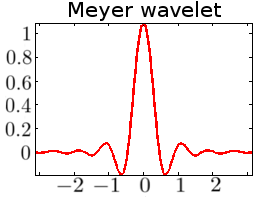}
    \caption{An example of mother wavelet, the Meyer wavelet.}
    \label{fig: meyer}
\end{center}\end{figure}

A parameterization that allows local refinements and local dropouts is the wavelet.
Wavelet is the basis developed for multi-resolution analysis (MRA) \cite{wavelet mallat}.
MRA is an increasing sequence of closed function spaces $\{V_j\}_{j\in \mathbb{Z}}$,
$$\cdots \subset V_{-1} \subset V_0 \subset V_1 \subset \cdots \,,$$
which can approximate functions with increasing resolutions as $j$ increases.
For univariate MRA, $V_j$'s satisfy the following properties known as
self-similarity:
\begin{equation*}\begin{split}
    &f(u) \in V_j \Leftrightarrow f(2u) \in V_{j+1}, \; j\in \mathbb{Z}\\
    &f(u) \in V_j \Leftrightarrow f(u-\frac{\eta}{2^j}) \in V_{j},
                    \; j\in \mathbb{Z},\, \eta\in \mathbb{Z}
\end{split}\end{equation*}
The wavelet bases for $V_j$ are given by
\begin{equation}
    \hat{\phi}_{j,\eta}(u) = 2^{j/2} \hat{\phi}(2^j u-\eta) \,,\quad \eta \in \mathbb{Z}
    \label{eqn: self similar wavelet}
\end{equation}
where $\hat{\phi}$ is called the mother wavelet satisfying 
$\hat{\phi}(u) \rightarrow 0$ for $ u \rightarrow -\infty $ and $\infty$.
An example mother wavelet, the Meyer wavelet, is shown in Figure \ref{fig: meyer}.

Because the gray-box model is unchanged by adding a constant to ${F}$,
the derivative of $\tilde{F}$, instead of $\tilde{F}$ itself, should be approximated.
Thereby, the bases for $\tilde{F}$ shall be chosen as the integrals of the wavelets, i.e.
\begin{equation}
    \phi_{j,\eta}(u) = \int_{-\infty}^u \hat{\phi}_{j,\eta}(u^\prime) du^\prime\,.
    \label{eqn: integral wavelet}
\end{equation}
$\phi_{j,\eta}$'s are sigmoid functions which satisfy
\begin{equation}
    \phi_{j,\eta}(u) = \left\{
        \begin{split}
            0,&\; u\rightarrow -\infty\\
            1,&\; u\rightarrow \infty
        \end{split}\right.
\end{equation}
It's easy to show $\phi_{j,\eta}$ also satisfies self-similarity, and can be 
represented by a mother sigmoid $\phi$,
\begin{equation}
    {\phi}_{j,\eta}(u) = {\phi}(2^j u-\eta) \,,\quad j\in \mathbb{Z} \,,\;\eta \in \mathbb{Z}
    \label{eqn: self similar sigmoid}
\end{equation}

\begin{figure}\begin{center}
    \includegraphics[width=4cm, height=3cm]{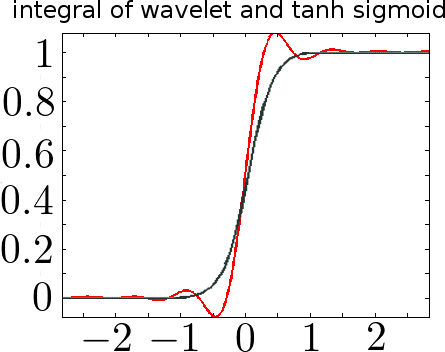}
    \caption{Red line: the integral \eqref{eqn: integral wavelet} of the Meyer wavelet.
             Black line: the logistic sigmoid function.}
    \label{fig: sigmoid}
\end{center}\end{figure}

There are many types of sigmoid functions.
The article will use the logistic sigmoid function as the mother sigmoid,
\begin{equation}
    \phi(u) = \frac{1}{1+ e^{-u}}\,.
    \label{eqn: logistic sigmoid}
\end{equation}

If $\tilde{F}$ is univariate, the logistic sigmoids $\phi_{j,\eta}$'s are used as the bases.
If $\tilde{F}$ is multivariate, the basis can be formed by the tensor product 
of univariate sigmoids,
\begin{equation}
    \phi_{\boldsymbol{j}, \boldsymbol{\eta}} (u_1, \cdots, u_k) = \phi_{j_1, \eta_1}(u_1)\cdots
    \phi_{j_k, \eta_k}(u_k)\,,
    \label{eqn: tensor basis}
\end{equation}
where $\boldsymbol{j}=(j_1, \cdots, j_k)\in \mathbb{Z}^k$, 
$\boldsymbol{\eta} = (\eta_1, \cdots, \eta_k) \in \mathbb{Z}^k$. To sum up,
\begin{equation}
    \tilde{F}(\cdot) = \sum_{\boldsymbol{j}, \boldsymbol{\eta}} 
                       \alpha_{\boldsymbol{j}, \boldsymbol{\eta}}
                       \phi_{\boldsymbol{j}, \boldsymbol{\eta}} (\cdot)\,,
    \label{eqn: linear expansion}
\end{equation}
where $\alpha$'s are the coefficients of the bases. A systematic procedure for choosing
a set of $\{\boldsymbol{j}, \boldsymbol{\eta}\}$
is presented in Section \ref{sec: adaptive basis}.
\\

\begin{figure}
    \begin{center}
        \includegraphics[width=5cm]{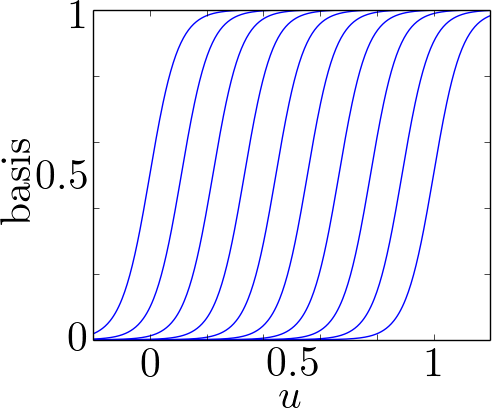}
        \caption{An ad hoc set of bases}
        \label{fig: sigmoid basis ad hoc}
    \end{center}
\end{figure}

For illustration, consider a numerical example where the gray-box model solves the 1-D Buckley-Leverett equation
\cite{Buckley Leverett}
\begin{equation}
    \frac{\partial u}{\partial t} + \frac{\partial}{\partial x}\Big(\underbrace{
    \frac{u^2}{1+ 2(1-u)^2}}_{F} \Big) = c\,,
    \label{eqn: Buckley-Leverett}
\end{equation}
with the initial condition $u(0,x)=u_0(x)$ and the periodic boundary condition $u(t,0)=u(t,1)$. 
Let $0 \le u_0(x) \le 1$ for all $x\in [0,1]$. This is because the 
Buckley-Leverett equation models the two-phase porous media flow where $u$ stands for the
saturation of a phase, and the saturation is always positive and no larger than one.
$c\in \mathbb{R}$ is a constant-valued control. $F$ is assumed unknown and is inferred by a twin model.
The twin model solves 
\begin{equation}
    \frac{\partial \tilde{u}}{\partial t} + \frac{\partial}{\partial x}\tilde{F}(\tilde{u})
    = c\,,
    \label{eqn: Buckley-Leverett twin}
\end{equation}
with the same $c$ and the same initial and boundary conditions. $\tilde{F}$ is parameterized by
\eqref{eqn: linear expansion} where
$j$'s and $\eta$'s are chosen ad hoc.
Figure \ref{fig: sigmoid basis ad hoc} gives an example of the bases used in this section.
To ensure the well-posedness of \eqref{eqn: twin min mismatch}, an 
ad hoc $L_1$ regularization on $\alpha$ is used in minimizing $\mathcal{M}$.\\

\begin{figure}\begin{center}
    \begin{subfigure}[t]{.4\textwidth}
        \centering
        \includegraphics[width=5cm]{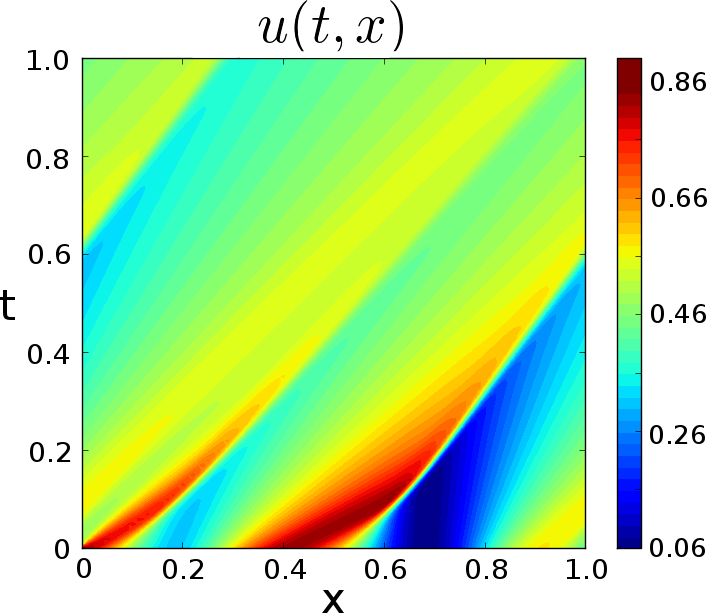}
        \caption{Gray-box model.}
        \label{fig: leftcol}
    \end{subfigure}
    \begin{subfigure}[t]{.4\textwidth}
        \centering
        \includegraphics[width=5cm]{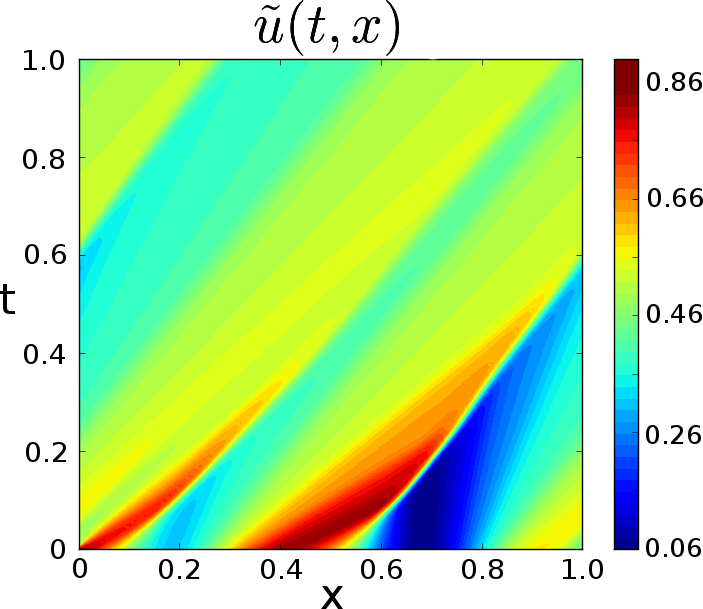}
        \caption{Trained twin model.}
        \label{fig: rightcol}
    \end{subfigure}
    \caption{Space time solutions.}
\end{center}\end{figure}

Figure \ref{fig: leftcol} shows the discretized space-time solution of \eqref{eqn: Buckley-Leverett}
for $x\in[0,1], \, t\in[0,1]$ and $c=0$. 
The solution is used to train a twin model according to \eqref{eqn: twin min mismatch}.
The discretized space-time solution of the trained twin model is shown in 
Figure \ref{fig: rightcol}.\\

\begin{figure}\begin{center}
    \includegraphics[width=6.5cm]{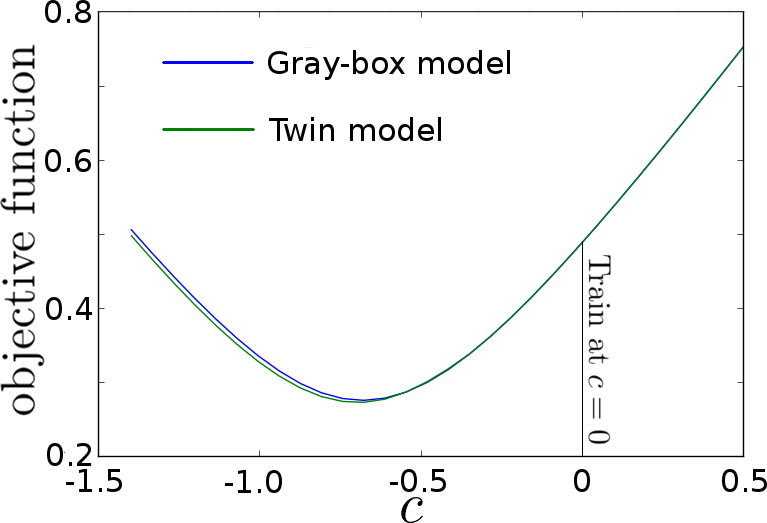}
    \caption{The objective function $\xi$ evaluated by either the gray-box model
             or the trained twin model.}
    \label{fig: objective ad hoc}
\end{center}\end{figure}

Once a twin model is trained, its adjoint can be used for gradient estimation.
Consider an objective function 
\begin{equation}
    \xi(c) \equiv \int_{x=0}^1 \big(u(1,x; c) - \frac{1}{2}\big)^2 \,\textrm{d}x\,.
    \label{eqn: objective ad hoc}
\end{equation}
Its gradient $\frac{d\xi}{dc}$ can be estimated by the trained twin model.
Figure \ref{fig: objective ad hoc} shows the objective function, evaluated using 
the gray-box model and the trained twin model. It is observed that the gradients of $\xi$
match closely at $c=0$, i.e. the control where the twin model is trained.\\

\begin{figure}
\begin{center}
    \includegraphics[width=13cm]{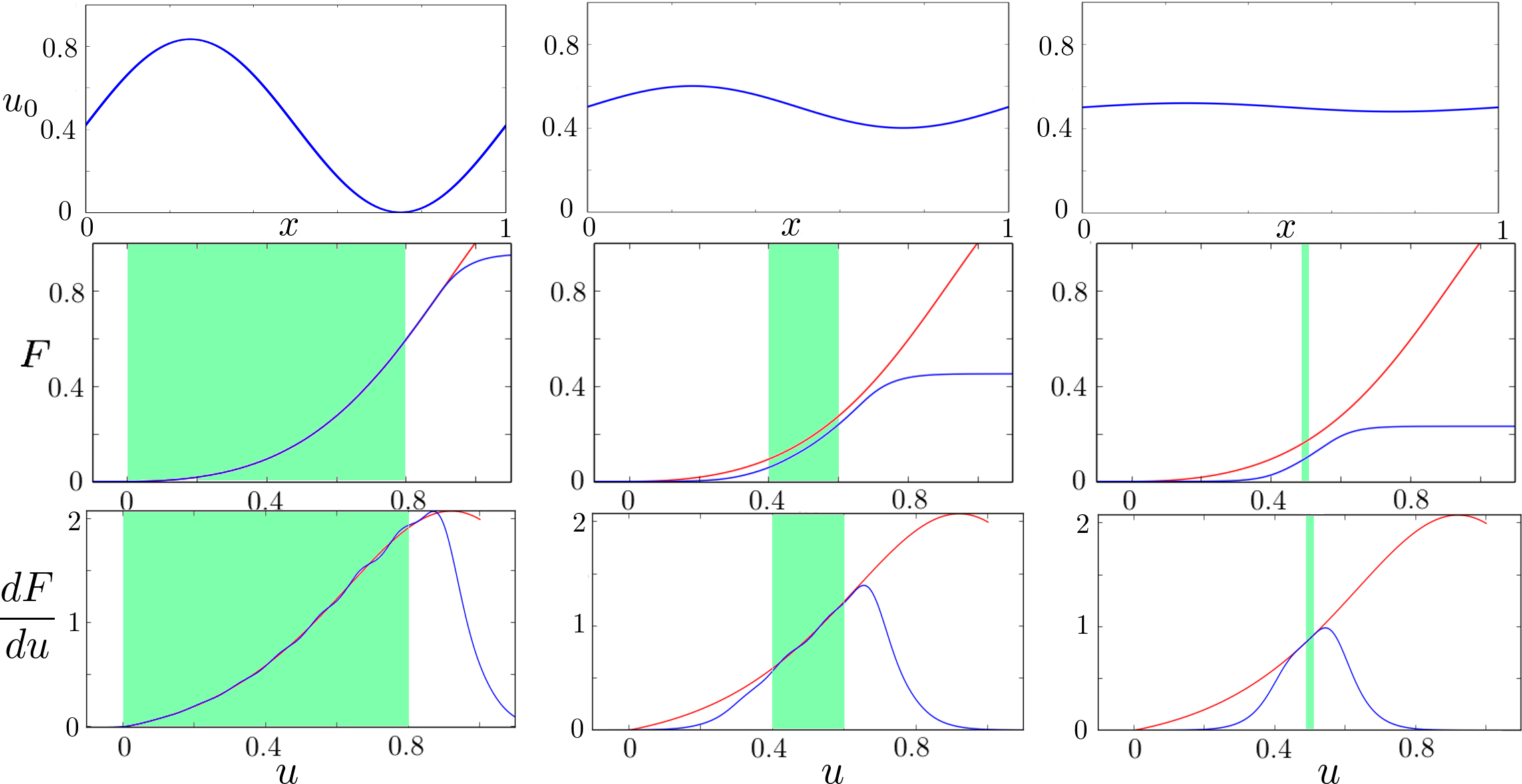}
    \caption{The first row shows the three different initial conditions used to generate
             the gray-box space-time solution. The second row compares the trained $\tilde{F}$ 
             (blue) and the Buckley-Leverett $F$ (red). The third row compares
             the trained $\frac{d\tilde{F}}{du}$ (blue) and the Buckley-Leverett 
             $\frac{dF}{du}$ (red). The green background highlights the domain of $u$ where
             the gray-box space-time solution exists.}
    \label{fig: combine 3}
\end{center}
\end{figure}

In addition to the gradient estimation, the inferred $\tilde{F}$ is examined.
If different solutions are used in the training, it is expected
that the trained twin model  will also be different.
Figure \ref{fig: combine 3} shows the training results for three different initial conditions.
Some observations can be made: 1) As expected, $\tilde{F}$ can differ from $F$ by a constant without
affecting $\mathcal{M}$; 2) $\frac{d\tilde{F}}{du}$ matches $\frac{dF}{du}$ only in a
domain of $u$ where the solution exists (indicated by the green area);
3) Sometimes the bases seem redundant thus can be safely dropped out. The issue seems particularly 
important in the third column, where most bases are suppressed;
4) Sometime the bases seem too coarse thus may be refined in order to reduce the minimal $\mathcal{M}$. 
The issue seems particularly important in the first column, where $\frac{d\tilde{F}}{du}$ exhibits
a wavy deviation from $\frac{dF}{du}$. 
Addressing these issues systematically is crucial to the rigorous development of
the twin model method.\\

\subsection{Elements for Adaptive Basis Construction}
\label{sec: adaptive basis}
This section develops several key 
elements that lead to the adaptive basis construction for twin models.
The heuristics for the adaptive basis construction, discussed in Section
\ref{sec: adaptive basis review}, are applied to build up a basis dictionary consisted of
only the significant candidates. The section is organized as follows:
Firstly, a formulation is provided to efficiently assess the 
significance of each candidate basis; Secondly, 
the neighborhood of a sigmoid basis is defined;
Thirdly, a metric is devised that determines when to add or remove a candidate basis.
The three elements are then employed to build the twin model algorithm in Section \ref{sec: twin algo}.\\

Given a basis dictionary $\boldsymbol{\phi}_\mathcal{A} = 
\{\phi_i\}_{i\in \mathcal{A}}$, define the ``minimal mismatch''
\begin{equation}
    \mathcal{M}^*(\mathcal{A}) = \min_{\boldsymbol{\alpha}_\mathcal{A} \in \mathbb{R}^{|\mathcal{A}|}}
    \mathcal{M}\left( \sum_{i\in \mathcal{A}}
    {\alpha}_i {\phi}_i \right)\,,
\end{equation}
to be the minimal solution mismatch \eqref{eqn: solution mismatch} if $\tilde{F}$ were parameterized
by $\boldsymbol{\phi}_\mathcal{A}$. $\mathcal{A}$ is a set containing $\{\boldsymbol{j},
\boldsymbol{\eta}\}$'s.
$\boldsymbol{\alpha}_{\mathcal{A}} = \{\alpha_i\}_{i\in \mathcal{A}}$ is the coefficient for 
$\boldsymbol{\phi}_\mathcal{A}$. 
Let $\boldsymbol{\alpha}_{\mathcal{A}}^* =\{\alpha_i^*\}_{i\in \mathcal{A}}$
be the optimal coefficients, and let
$\tilde{F}^*_\mathcal{A} = \sum_{i\in \mathcal{A}} \alpha^*_i \phi_i$.
Consider appending $\boldsymbol{\phi}_{\mathcal{A}}$ 
by an additional basis
$\phi_l$, and let $\boldsymbol{\phi}_{\mathcal{A}^\prime} = \left\{
\boldsymbol{\phi}_{\mathcal{A}}, \phi_l \right\}$,
$\mathcal{A}^\prime = \left\{ \mathcal{A}, l \right\}$.
The minimal mismatch for the appended basis dictionary
$\boldsymbol{\phi}_{\mathcal{A}^\prime}$ is
\begin{equation}
    \mathcal{M}^*(\mathcal{A}^\prime) 
    = \min_{\boldsymbol{\alpha}_{\mathcal{A}^\prime} \in \mathbb{R}^{|\mathcal{A}|+1}}
    \mathcal{M}\left( \sum_{i\in \mathcal{A}^\prime}
    {\alpha}_{i} {\phi}_{i} \right)\,,
\end{equation}
Clearly $\mathcal{M}^*(\mathcal{A}^\prime) \le \mathcal{M}^*(\mathcal{A})$.
Define the ``mismatch improvement'' to be
\begin{equation}
    \Delta \mathcal{M}^*\left(\mathcal{A}, l\right) = \mathcal{M}^*(\mathcal{A}) - 
    \mathcal{M}^*(\mathcal{A}^\prime)
    \label{eqn: mismatch improvement}
\end{equation} 
Approximate \eqref{eqn: mismatch improvement} by Taylor expansion, we get
\begin{equation}
    \Delta\mathcal{M}^*\left(\mathcal{A}, l\right) 
    \approx -\left(\int_{u\in \mathbb{R}^k} \left.\frac{d\mathcal{M}}{d \tilde{F}}
    \right|_{\tilde{F}_\mathcal{A}^*} \phi_l
    \; \textrm{d} u \right) \alpha_l \,,
    \label{eqn: taylor expansion}
\end{equation}
where $\frac{d\mathcal{M}}{d\tilde{F}}$ is the derivative of $\mathcal{M}(\tilde{F})$ 
with respect to $\tilde{F}$, evaluated on $\tilde{F} = \tilde{F}^*_\mathcal{A}$.
For a twin model consisted of a system of $k$ equations,
$\tilde{F}$ is a function of $u\in\mathbb{R}^k$, thus 
$\frac{d\mathcal{M}}{d\tilde{F}}$ is also a function of $u\in\mathbb{R}^k$.
As discussed in the previous sections, $\frac{d \mathcal{M}}{d\tilde{F}}$ is non-zero only in a
domain where there is solution. Thus \eqref{eqn: taylor expansion} 
can be integrated by quadrature only over a bounded domain.
The absolute value of the coefficient for $\alpha_l$,
\begin{equation}
    s_l(\mathcal{A}) \equiv \left|\int_{u\in \mathbb{R}^k} \left.\frac{d\mathcal{M}}{d \tilde{F}}
    \right|_{\tilde{F}_\mathcal{A}^*} \phi_l \; \textrm{d} u \right|\,,
    \label{eqn: basis significance}
\end{equation}
estimates the significance of the basis $\phi_l$ \cite{weight selection}. 
If there are multiple candidate bases, \eqref{eqn: basis significance} can be used to rank their 
significance.\\

\begin{figure}\begin{center}
    \begin{subfigure}[t]{.48\textwidth}
        \centering
        \includegraphics[width=7.5cm]{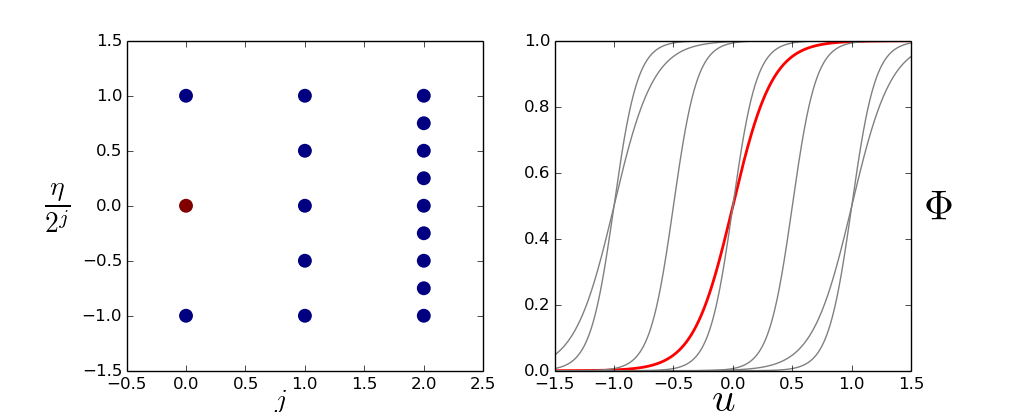}
        \caption{$\left(0,\frac{0}{2^0}\right)$}
        \label{fig: basis 0}
    \end{subfigure}
    \begin{subfigure}[t]{.48\textwidth}
        \centering
        \includegraphics[width=7.5cm]{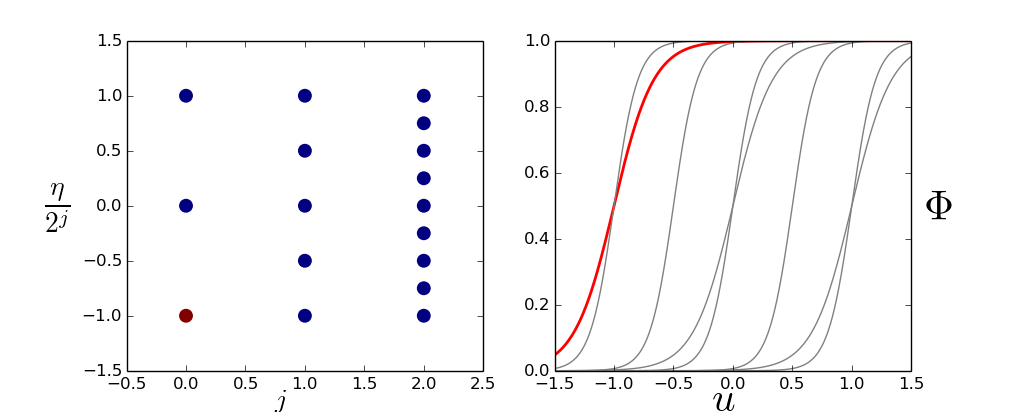}
        \caption{$\left(0, \frac{-1}{2^0}\right)$}
        \label{fig: basis 1}
    \end{subfigure}
    \begin{subfigure}[t]{.48\textwidth}
        \centering
        \includegraphics[width=7.5cm]{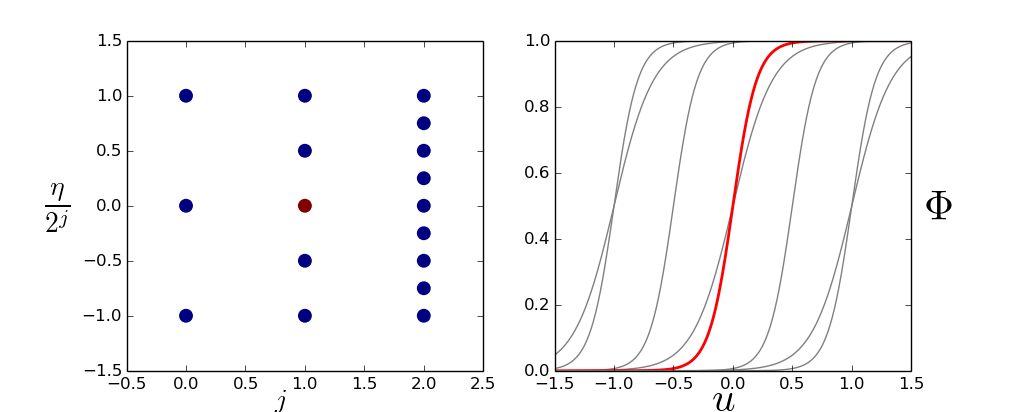}
        \caption{$\left(1, \frac{1}{2^0}\right)$}
        \label{fig: basis 2}
    \end{subfigure}
    \begin{subfigure}[t]{.48\textwidth}
        \centering
        \includegraphics[width=7.5cm]{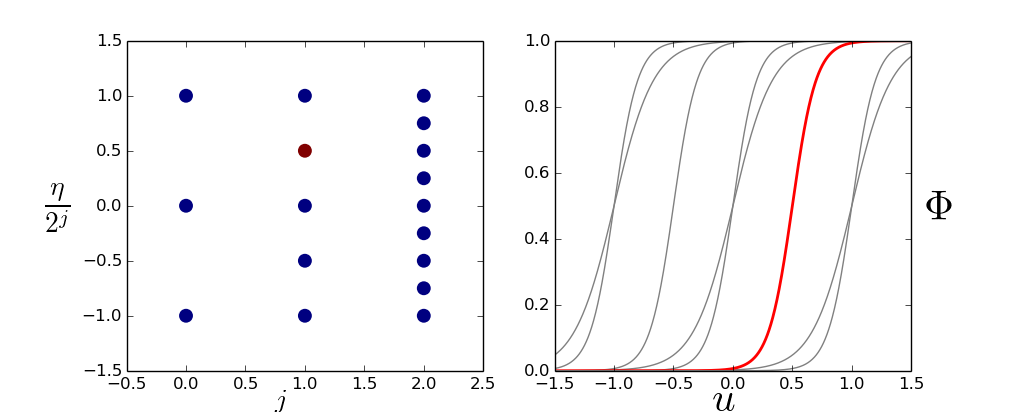}
        \caption{$\left(1, \frac{1}{2^1}\right)$}
        \label{fig: basis 3}
    \end{subfigure}
    \caption{An illustration of the tuple representation and the corresponding univariate sigmoid.}
\end{center}\end{figure}

In the sequel, a compact representation of the sigmoid bases is introduced.
The univariate basis function, $\phi_{j,\eta}$ in
\eqref{eqn: self similar sigmoid}, is represented by a tuple $(j, \frac{\eta}{2^j})$,
where $j$ can be viewed as the ``resolution'', and $\frac{\eta}{2^j}$ is
the center of the basis.
Similarly, the $k$-variate basis function, $\phi_{\boldsymbol{j}, \boldsymbol{\eta}}$ in
\eqref{eqn: tensor basis}, is represented by a tuple $\left(\boldsymbol{j}, \,
\frac{\boldsymbol{\eta}}{2^{\boldsymbol{j}}}\right) 
= \left(\{j_1,\cdots, j_k\}, \, \left\{\frac{\eta_1}{2^{j_1}}, \cdots, \frac{\eta_k}{2^{j_k}}\right\}\right)$.
Thus, a sigmoid can be represented by a point in a $2k$-dimensional space.
The representation is illustrated in Figure \ref{fig: basis 0} thru. \ref{fig: basis 3} 
for the univariate case.

\begin{figure}\begin{center}
    \begin{subfigure}[p]{1.\textwidth}
        \centering
        \includegraphics[width=10cm]{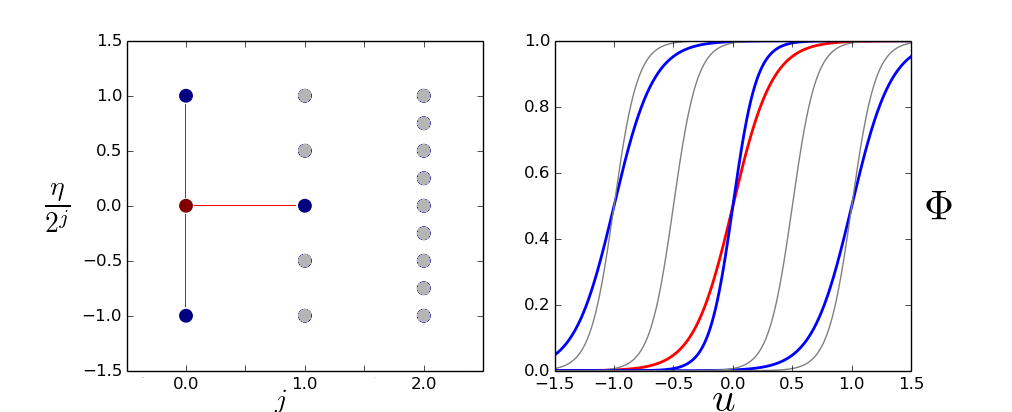}
        \caption{$\mathcal{N}\left[\left(0,\frac{0}{2^0}\right)\right]$}
        \label{fig: basis neighbor}
    \end{subfigure}
    \begin{subfigure}[p]{1.\textwidth}
        \centering
        \includegraphics[width=10cm]{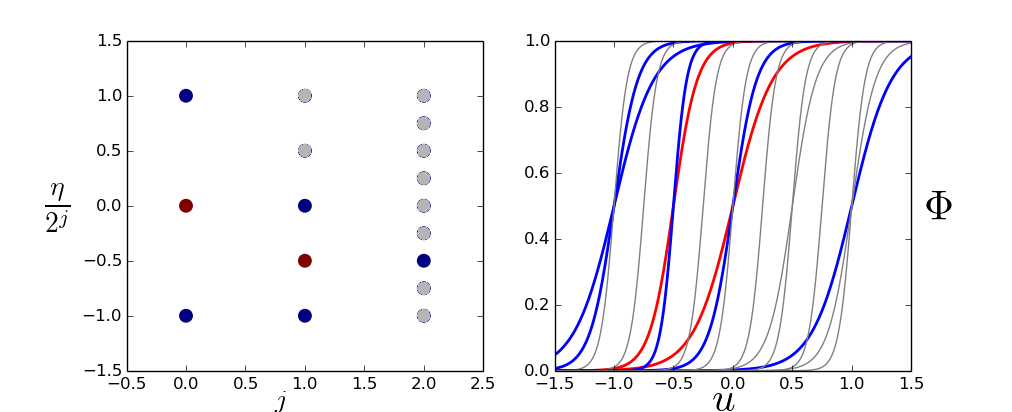}
        \caption{$\mathcal{N}\left[ \big(0, \frac{0}{2^0}\big), \big(1, \frac{-1}{2^1}\big)
                 \right]$}
        \label{fig: union neighbor}
    \end{subfigure}
    \caption{Neighborhood for univariate bases. 
             $(a)$ shows the neighborhood (blue)
             of a single basis (red).  $(b)$ shows the neighborhood (blue) 
             of several bases (red). 
             The left column represents the basis on the $\left(j, \frac{\eta}{2^j}\right)$ plane,
             and the right column shows the actual basis $\phi_{j,\eta}$.}
\end{center}\end{figure}

Using this representation, define the ``neighborhood'' of a univariate sigmoid 
$(j,\frac{\eta}{2^j})$ to be 
\begin{equation}
    \mathcal{N}\left[ \left(j,\frac{\eta}{2^j}\right) \right]
    = \left\{ 
        \left( j+1, \frac{\eta}{2^j} \right),\,
        \left( j, \frac{\eta\pm 1}{2^j} \right)
    \right\}\,.
    \label{eqn: neighborhood 1D}
\end{equation}
The neighborhood contains 1) a basis 
$\left( j+1, \frac{\eta}{2^j} \right)$ 
with an increment of resolution; and
2) two basis $\left( j, \frac{\eta\pm 1}{2^j} \right)$
with the same resolution but a marginal shift of center.
For illustration, the neighborhood of $\left(0,\frac{0}{2^0}\right)$
is shown in Figure \ref{fig: basis neighbor}.
Similarly, define the neighborhood of a multivariate sigmoid to be
\begin{equation}\begin{split}
    &\mathcal{N}
    \left[
         \left(
               \boldsymbol{j}, \frac{\boldsymbol{\eta}}{2^{\boldsymbol{j}}}
         \right)
    \right] = 
    \mathcal{N}\left[ \left(\{j_1, \cdots, j_k\} , \left\{
    \frac{\eta_1}{2^{j_1}}, \cdots, \frac{\eta_k}{2^{j_k}} \right\} \right) \right]\\
    = & \bigg\{
            \left( \left\{ j_1+1,\cdots, j_k\right\},
                   \left\{ \frac{\eta_1}{2^{j_1+1}}, \cdots, \frac{\eta_k}{2^{j_k}} \right\}
            \right) \, \cdots \, ,
            \left( \left\{ j_1,\cdots, j_k+1\right\},
                   \left\{ \frac{\eta_1}{2^{j_1}}, \cdots, \frac{\eta_k}{2^{j_k+1}} \right\}
            \right),\,   \\
          & 
            \left( \left\{ j_1,\cdots, j_k\right\},
                   \left\{ \frac{\eta_1\pm 1}{2^{j_1}}, \cdots, \frac{\eta_k}{2^{j_k}} \right\}
            \right) \, \cdots \, ,
            \left( \left\{ j_1,\cdots, j_k\right\},
                   \left\{ \frac{\eta_1}{2^{j_1}}, \cdots, \frac{\eta_k\pm 1}{2^{j_k}} \right\}
            \right) \bigg\}\,,
    \label{eqn: neighborhood kD}
\end{split}\end{equation}
which consists of $k$ bases with incremental resolution, and $2k$ bases with center shifts.
It is easy to see that a basis $(\boldsymbol{j}_0, \frac{\boldsymbol{\eta}_0}{2^{\boldsymbol{j}_0}})$ 
can be connected to any basis $(\boldsymbol{j}, \frac{\boldsymbol{\eta}}{2^{\boldsymbol{j}}})$ 
with $\boldsymbol{j} \ge \boldsymbol{j}_0$ through a chain
of neighborhoods.
In addition, define the neighborhood of multiple sigmoid functions to be the union
of each individual's neighborhood, as illustrated by Figure \ref{fig: union neighbor}.
\begin{equation}
    \mathcal{N}\left[(\boldsymbol{j}_1, \frac{\boldsymbol{\eta}_1}{2^{\boldsymbol{j}_1}} ), \cdots, 
    ( \boldsymbol{j}_n, \frac{\boldsymbol{\eta}_n}{2^{\boldsymbol{j}_n}}) \right]
    = \mathcal{N}\left[(\boldsymbol{j}_1, \frac{\boldsymbol{\eta}_1}{2^{\boldsymbol{j}_1}})\right]\bigcup \cdots 
      \bigcup \mathcal{N}\left[(\boldsymbol{j}_n, \frac{\boldsymbol{\eta}_n}{2^{\boldsymbol{j}_n}})\right]\,.
\end{equation}\\

\begin{figure}
    \begin{center}
        \includegraphics[width=3.2cm]{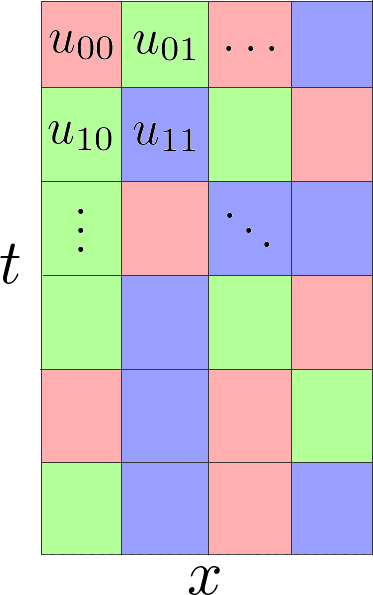}
        \caption{The discretized 
                 gray-box solution is shuffled into $3$ sets, each indicated by a color.}
        \label{fig: shuffle}
    \end{center}
\end{figure}

Although the mismatch improvement, $\Delta \mathcal{M}^*\left(\mathcal{A}, l\right)$, is always
non-negative, it is inadvisable to cram the basis dictionary with too many bases,
otherwise a twin model can be overfitted.
Therefore, a criterion is required to determine if a candidate basis shall be added to or 
removed from the basis dictionary. This can be achieved by cross validation,
in particular, $k$-fold cross validation \cite{cross validation}.
Given a basis dictionary, the $k$-fold cross validation proceeds in the following three steps:
Firstly, the gray-box solution $\boldsymbol{u}$ is shuffled randomly into $k$ disjoint sets
$\left\{\boldsymbol{u}_1 , \boldsymbol{u}_2, \cdots, \boldsymbol{u}_k\right\}$.
An illustration for $k=2$ is shown in Figure \ref{fig: shuffle}.

Secondly, $k$ twin models who share the same basis dictionary 
are trained so their space-time solutions match
all but one sets, shown in \eqref{eqn: cross validation}.
$T_i$ indicates the $i$th twin model.
\begin{equation}\begin{split}
T_1 &= \texttt{TrainTwinModel}(\boldsymbol{u}_2, \boldsymbol{u}_3, \cdots, \boldsymbol{u}_k)\\
T_2 &= \texttt{TrainTwinModel}(\boldsymbol{u}_1, \boldsymbol{u}_3, \cdots, \boldsymbol{u}_k)\\
&\cdots\\
T_k &= \texttt{TrainTwinModel}(\boldsymbol{u}_1, \boldsymbol{u}_2, \cdots, \boldsymbol{u}_{k-1})
\label{eqn: cross validation}
\end{split}\end{equation}
Thirdly, each trained twin model is validated on the remaining set. In particular,
the solution mismatch for the validation set is computed, as shown in 
\eqref{eqn: cross validation mismatch}.
\begin{equation}\begin{split}
    \mathcal{M}_1 &= \texttt{MismatchValidation}\left( T_1 , \boldsymbol{u}_1\right)\\
    \mathcal{M}_2 &= \texttt{MismatchValidation}\left( T_2 , \boldsymbol{u}_2\right)\\
    &\cdots\\
    \mathcal{M}_k &= \texttt{MismatchValidation}\left( T_k , \boldsymbol{u}_k\right)
    \label{eqn: cross validation mismatch}
\end{split}\end{equation}
The mean value of validation errors,
\begin{equation}
    \overline{\mathcal{M}} = \frac{1}{k}\left(\mathcal{M}_1 + \mathcal{M}_2 + \cdots 
    + \mathcal{M}_k\right)
\end{equation}
measures the performance of the basis dictionary. A basis shall be added to or removed from
the dictionary only if such action reduces $\overline{\mathcal{M}}$. In practice, cross validation
proliferates the computational cost. 
Therefore, a small $k$ is preferrable if cost is a concern. 
All the numerical examples in this article use $k=2$.\\

\subsection{Algorithm}
\label{sec: twin algo}
Based upon the developments in the previous sections, a twin model algorithm 
with adaptive basis construction is devised.\\

\begin{algorithm}[H]
\begin{algorithmic}[1]
\REQUIRE{Initial basis dictionary $\boldsymbol{\phi}_\mathcal{A}$, 
        coefficients $\alpha_{\mathcal{A}}= \mathbf{0}$},
        Validation error $\overline{\mathcal{M}}_0 = \infty$,
        Gray-box solution $\boldsymbol{u}$.
\STATE Minimize solution mismatch 
       $
           \alpha_{\mathcal{A}} \leftarrow \argmin_{\alpha} \mathcal{M}\left(\sum_{i \in \mathcal{A}} \alpha_i \phi_i\right)
       $ 
\LOOP 
\STATE Find $\phi_{l}\in \mathcal{N}(\boldsymbol{\phi}_{\mathcal{A}}) \backslash 
       \boldsymbol{\phi}_{\mathcal{A}}$ with the maximal $s_l(\mathcal{A})$ \\
       $\mathcal{A} \leftarrow 
       \mathcal{A} \bigcup \{l\}$, $\boldsymbol{\phi}_{\mathcal{A}} \leftarrow 
       \boldsymbol{\phi}_{\mathcal{A}}
       \bigcup \{\phi_{l}\}$, $\alpha_l = 0$,
       $\alpha_{\mathcal{A}} \leftarrow \{\alpha_{\mathcal{A}}, \alpha_l\}$
\STATE Compute $\overline{\mathcal{M}}$ by $k$-fold cross validation.
\IF {$\overline{\mathcal{M}} < \overline{\mathcal{M}}_0$} 
    \STATE $\overline{\mathcal{M}}_0\leftarrow \overline{\mathcal{M}}$\\
       $
           \alpha_{\mathcal{A}} \leftarrow \argmin_{\alpha} \mathcal{M}\left(\sum_{i \in \mathcal{A}} \alpha_i \phi_i\right)
       $
\ELSE \STATE 
       $\mathcal{A} \leftarrow 
       \mathcal{A} \backslash \{l\}$, $\boldsymbol{\phi}_{\mathcal{A}} \leftarrow 
       \boldsymbol{\phi}_{\mathcal{A}}
       \backslash \{\phi_{l}\}$, $\alpha_{\mathcal{A}} \leftarrow \alpha_{\mathcal{A}} \backslash 
       \{\alpha_l\}$
      \textbf{break}
\ENDIF
\STATE Find $\phi_{l^\prime}\in 
       \boldsymbol{\phi}_{\mathcal{A}}$ with the least $s_{l^\prime}(\mathcal{A})$  \\
\IF {${l^\prime} \neq {l}$}
\STATE  
       $\mathcal{A} \leftarrow 
       \mathcal{A} \backslash \{l^\prime\}$,
       $\boldsymbol{\phi}_{\mathcal{A}} \leftarrow \boldsymbol{\phi}_{\mathcal{A}}
        \backslash\{\phi_{l^\prime}\}$, 
       $\alpha_{\mathcal{A}} \leftarrow \alpha_{\mathcal{A}}\backslash \{\alpha_{l^\prime}\}$
\STATE Compute $\overline{\mathcal{M}}$ by $k$-fold cross validation.
\IF{  $\overline{\mathcal{M}} < \overline{\mathcal{M}}_0$ }
\STATE $\overline{\mathcal{M}}_0 \leftarrow \overline{\mathcal{M}}$\\
       $
           \alpha_{\mathcal{A}} \leftarrow \argmin_{\alpha} \mathcal{M}\left(\sum_{i \in \mathcal{A}} \alpha_i \phi_i\right)
       $
\ELSE \STATE
       $\mathcal{A} \leftarrow \mathcal{A}\bigcup \{l^\prime\}$,
       $\phi_{\mathcal{A}}\leftarrow \phi_{\mathcal{A}} 
       \bigcup \{\phi_{l^\prime}\}$, $\alpha_{\mathcal{A}} \leftarrow 
       \alpha_{\mathcal{A}}\bigcup \{\alpha_{l^\prime}\}$
\ENDIF
\ENDIF
\ENDLOOP
\ENSURE $\mathcal{A}$, $\phi_{\mathcal{A}}$, $\alpha_{\mathcal{A}}$.
\end{algorithmic}
\caption{Training twin model with adaptive basis construction.}
\label{alg: train twin}
\end{algorithm}

Algorithm \ref{alg: train twin} adopts the heuristics of the forward-backward iteration
discussed in Section \ref{sec: adaptive basis review}.
The algorithm starts from training a twin model using 
a simple basis dictionary. Usually the starting
dictionary contains one basis for each dimension with very low resolution.
Details of the choice are given in Section 
\ref{sec: twin numerical results} along with numerical examples.
The main part of the algorithm iterates over a forward step (line 3-9) and a backward step (line
10-19). The forward step firstly finds the most promising candidate in the neighborhood of the current
dictionary for addition, according to \eqref{eqn: basis significance}.
If the addition indeed reduces the cross validation error, the candidate is appended
to the dictionary; otherwise it is rejected. 
If the basis is appended, the coefficients are updated by minimizing the solution mismatch, which can be
implemented by the Broyden-Fletcher-Goldfarb-Shannon (BFGS) algorithm \cite{quasiNewton}.
The backward step finds the most promising
candidate in the current dictionary for deletion. 
If the deletion reduces the cross validation error, the candidate
is removed from the dictionary. 
If the basis is deleted, the coefficients are updated by BFGS again.
The iteration exits when the most promising addition no longer 
reduces the validation error. In the end, the algorithm provides the basis dictionary and its
coefficients as the output.\\

The algorithm requires to train multiple twin models at each iteration. For
$k=2$, $6$ twin models are trained if both the forward and the backward step are acceptive.
In practice, the trained coefficients at the last iteration usually provide good initial guess
for the next iteration. Nonetheless, the algorithm can be costly if the dictionary turns out
to have a high cardinality which results in a large number of iterations before the dictionary 
construction completes. Therefore, a numerical shortcut is provided in Section 
\ref{sec: trunc error} that significantly reduces the cost.\\

\subsection{Minimizing the Truncation Error}
\label{sec: trunc error}
In the previous sections, a twin model is trained to minimize the solution mismatch.
The training can be expensive. Because the minimization
of the solution mismatch, coupled with the adaptive basis construction, can require
a large number of solution mismatch evaluations, and each evaluation involves one
twin model simulation. To reduce the computational cost, a ``pre-training'' step is proposed
where an ``integrated truncation error'' is minimized. 
A pre-trained twin model is then ``fine tuned'' to minimize the solution mismatch.
The applicability of the pre-training is studied; in particular, I study
under what condition can the solution mismatch be bounded by the integrated truncation error.
Finally, a stochastic gradient descent approach is adopted that efficiently 
minimizes the integrated truncation error.\\

Define
\begin{equation}
    \tau = \frac{\partial u}{\partial t} 
    + \nabla \cdot \big(D \tilde{F}(u)\big) - q(u,c)\,,
    \label{eqn: residual}
\end{equation}
which is the residual if the gray-box PDE's solution is plugged in the twin-model PDE 
\eqref{eqn: govern twin model}. Let its discretization be $\boldsymbol{\tau}$.
For simplicity, assume the gray-box simulator and its twin model use the same space-time 
discretization. 
$\boldsymbol{\tau}$ can be obtained by plugging the discretized gray-box solution in the
twin model simulator.
Define the integrated truncation error to be
\begin{equation}
    \mathcal{T}(\tilde{F}) = \sum_{i=1}^M \sum_{j=1}^N w_{ij} \boldsymbol{\tau}_{ij}^2 \,,
    \label{eqn: truncation error}
\end{equation}
where $w_{ij}$ are the same quadrature weights as in \eqref{eqn: solution mismatch}. 
$i,j$ are the indices for time and space discretization as in the previous sections.
I propose to pre-train a twin model using Algorithm \ref{alg: train twin} with $\mathcal{M}$
replaced by $\mathcal{T}$. In other words, in the pre-training step, the coefficients are determined by
\begin{equation}
    \alpha_{\mathcal{A}} \leftarrow \argmin_{\alpha} \mathcal{T}\left(\sum_{i \in \mathcal{A}}
    \alpha_i \phi_i\right)\,.
    \label{eqn: minimize truncation error}
\end{equation}
Besides, the estimator for the significance of a candidate basis, $s_l(\mathcal{A})$, is replaced by 
\begin{equation}
    s_l^t(\mathcal{A}) \equiv \left|\int_{u\in \mathbb{R}^k} \left.\frac{d\mathcal{T}}{d \tilde{F}}
    \right|_{\tilde{F}_\mathcal{A}^*} \phi_l \; \textrm{d} u \right|\,.
    \label{eqn: basis significance 2}
\end{equation}
Finally, the validation error, $\overline{\mathcal{M}}$, is replaced by
\begin{equation}
    \overline{\mathcal{T}} = \frac{1}{k}\left( \mathcal{T}_1 + \mathcal{T}_2 + \cdots + 
    \mathcal{T}_k\right)\,,
    \label{eqn: pre train validation error}
\end{equation}
where
\begin{equation}
    \mathcal{T}_i = \texttt{IntegratedTruncationError}(T_i, \boldsymbol{u}_i)\,.
\end{equation}
for $i=1, \cdots, k$.
Using the pre-trained basis dictionary $\phi^t_{\mathcal{A}}$,
the twin model is then fine tuned by minimizing the 
solution mismatch, where $\alpha^t_{\mathcal{A}}$ is used as the initial guess and is adjusted 
according to \eqref{eqn: twin min mismatch}. 
For a simulation with implicit schemes, the residual and the integrated truncation error 
are cheaper to evaluate than the solution mismatch, thereby the benefit of the pre-train. \\

However, $\mathcal{M}$ may not be bounded by
$\mathcal{T}$. A sufficient condition under which the bound exists is provided
by Theorem \ref{theorem: 2}.\\

\begin{theorem}
    Consider a twin model simulator whose one-step time marching is
    \begin{equation}
        \mathcal{G}_i:\, \mathbb{R}^N\mapsto\mathbb{R}^N,\, \tilde{\boldsymbol{u}}_{i\cdot}\rightarrow 
        \tilde{\boldsymbol{u}}_{i+1\cdot}=\mathcal{G}_i
        \tilde{\boldsymbol{u}}_{i\cdot} \,,\quad i=1,\cdots, M-1\,.
    \end{equation}
    Assume the quadrature weights are time-independent, i.e.
    $w_{ij} = w_{j}$ for all $i,j$.
    If $\mathcal{G}_i$ satisfies
    \begin{equation}
        \|\mathcal{G}_ia-\mathcal{G}_ib\|^2_{W} \le \beta \|a-b\|^2_{W} \,,
        \label{eqn: contractive}
    \end{equation}
    with $\beta<1$,
    for any $a, b \in \mathbb{R}^N$ and for all $i$,
    then 
    \begin{equation}
        \mathcal{M} \le \frac{1}{1-\beta} \mathcal{T}\,,
    \end{equation}
    where
    \begin{equation}
        \|v\|^2_{W} \equiv v^T 
            \begin{pmatrix}
                {w_{1}} && \\
                & \ddots & \\
                && {w_{N}}
            \end{pmatrix} v
    \end{equation}
    for any $v\in \mathbb{R}^N$.
    \label{theorem: 2}
\end{theorem}
The proof is given in Appendix \ref{proof 2}. If the twin model is a contractive dynamical 
system \cite{contractive system}, as given by \eqref{eqn: contractive}, then the
solution mismatch can be bounded by the integrated truncation error. In contrast, 
the bound may not exist for non-contractive dynamical systems, for example for systems 
that exhibit bifurcation. It is a future work to further investigate the applicability of the pre-training
theoretically, in particular, to investigate the necessary and sufficient condition for the bound. 
This article will explore the usefulness of the pre-training by several numerical test cases.\\

Because the residual $\boldsymbol{\tau}$ can be evaluated explicitly given the gray-box
solution, 
the evaluation can be decoupled for different space-time grid points $\{i,j\}$. By viewing
the truncation error at each $\{i,j\}$ as a stochastic sample, \eqref{eqn: minimize truncation error}
can be solved by stochastic gradient descent, Algorithm \ref{alg: 2}.\\
\begin{algorithm}
\begin{algorithmic}[1]
    \REQUIRE $\alpha = \alpha_0$
    \FOR{$(i,j)=(1,1)$ \TO $(M,N)$ }
         \IF {not converged}
             \STATE $\alpha \leftarrow \alpha -\lambda w_{ij}\left.\frac{\partial }{\partial \alpha} 
             \boldsymbol{\tau_{ij}} \right.$
         \ELSE
             \STATE \textbf{break}
         \ENDIF
    \ENDFOR
    \ENSURE $\alpha$
\end{algorithmic}
\caption{Minimizing the integrated truncation error by stochastic gradient descent.}
\label{alg: 2}
\end{algorithm}

$\lambda>0$ is a tunable step size. $\lambda$ can tuned manually to 
increase convergence speed while avoiding divergence \cite{stochastic search}.
In practice, it is beneficial to compute the gradient
against more than one grid points (called a ``mini-batch'') at each iteration. This is because
the code can take advantage of vectorization libraries rather than computing the residual 
at each grid point separately.\\

\section{Numerical Results}
\label{sec: twin numerical results}
This section demonstrates the twin model on the estimation of the gradients for
several numerical examples.

\subsection{Buckley-Leverett Equation}
\label{sec: chap 2 BL}
Section \ref{sec: flux param} has applied a sigmoid parameterization to the gray-box model
governed by the Buckley-Leverett equation \eqref{eqn: Buckley-Leverett}.
In this section, the same problem is studied but using the adaptive basis construction developed 
in Section \ref{sec: twin algo} and \ref{sec: trunc error}. 
The initial dictionary, $\boldsymbol{\phi}_{\mathcal{A}}$, is selected to contain a single basis
$\left(0, \frac{0}{2^0}\right)$. Clearly the choice is not unique. 
As long as the
initial basis has a low resolution and is centered around $\left[u_{\min}, u_{\max}\right]$,
Algorithm \ref{alg: train twin} shall build the dictionary adaptively.\\

\begin{figure}\begin{center}
    \begin{subfigure}[t]{.32\textwidth}
        \centering
        \includegraphics[width=4.5cm]{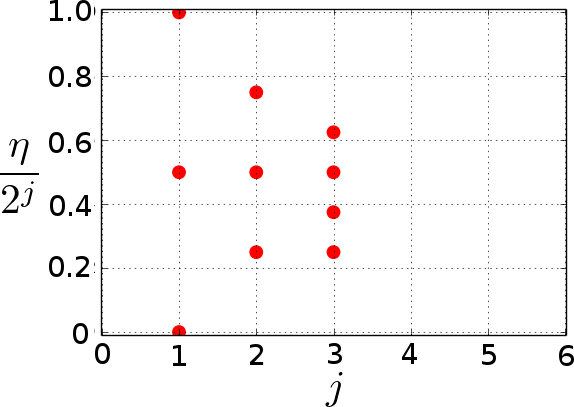}
        \caption{Solution 1}
        \label{fig: basis pnt 1}
    \end{subfigure}
    \begin{subfigure}[t]{.32\textwidth}     
        \centering
        \includegraphics[width=4.5cm]{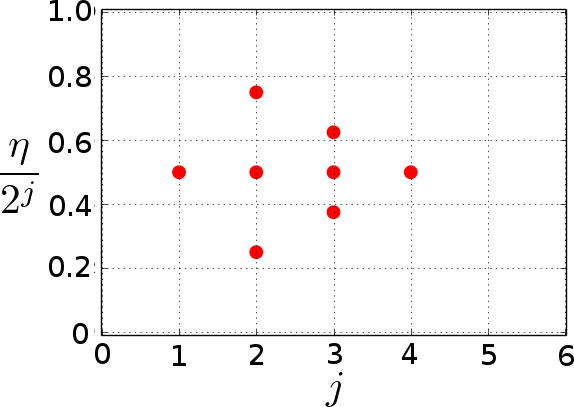}
        \caption{Solution 2}
        \label{fig: basis pnt 2}
    \end{subfigure}
    \begin{subfigure}[t]{.32\textwidth}
        \centering
        \includegraphics[width=4.5cm]{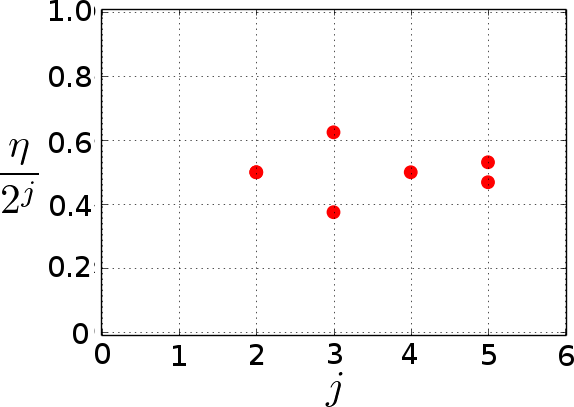}
        \caption{Solution 3}
        \label{fig: basis pnt 3}
    \end{subfigure}
    \caption{The basis dictionary for the three solutions in Figure \ref{fig: combine 3}.}
    \label{fig: basis pnt}
\end{center}\end{figure}

Figure \ref{fig: basis pnt} shows the selected bases
for the three solutions in Figure \ref{fig: combine 3}, respectively, by using
the pre-train step.
As $\left[u_{\min}, u_{\max}\right]$ shrinks, the dictionary's cardinality reduces and
the resolution increases. \\

\begin{figure}\begin{center}
    \begin{subfigure}[t]{.32\textwidth}
        \centering
        \includegraphics[width=4.8cm, height=4.1cm]{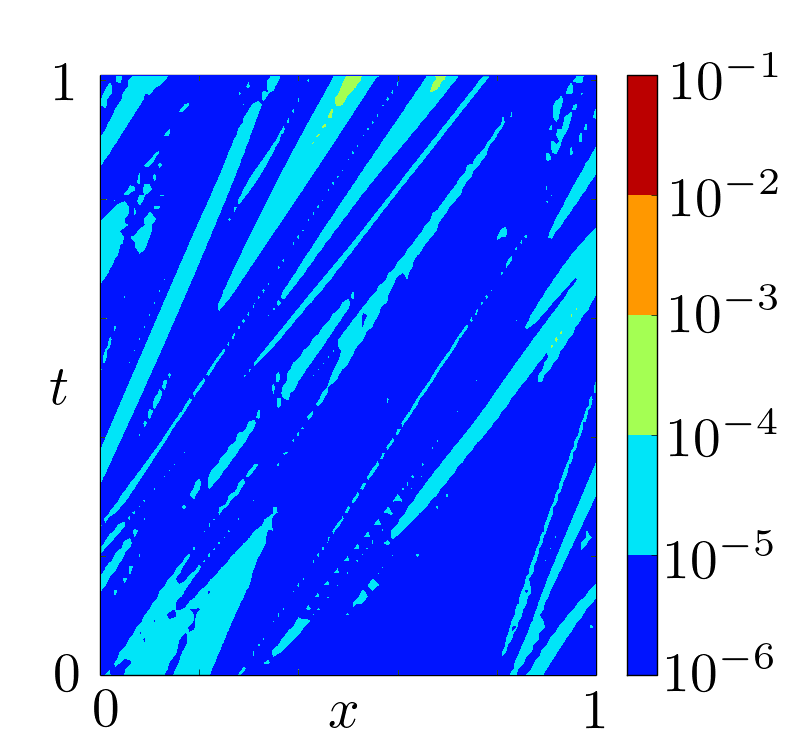}
        \caption{Solution 1}
        \label{fig: grad BL 1}
    \end{subfigure}
    \begin{subfigure}[t]{.32\textwidth}     
        \centering
        \includegraphics[width=4.8cm, height=4.1cm]{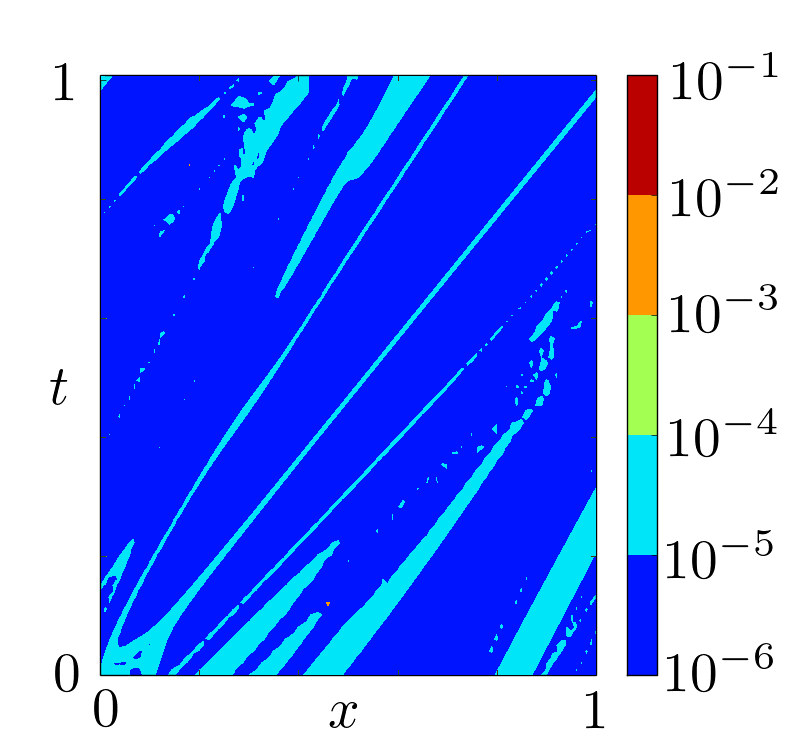}
        \caption{Solution 2}
        \label{fig: grad BL 2}
    \end{subfigure}
    \begin{subfigure}[t]{.32\textwidth}
        \centering
        \includegraphics[width=4.8cm, height=4.1cm]{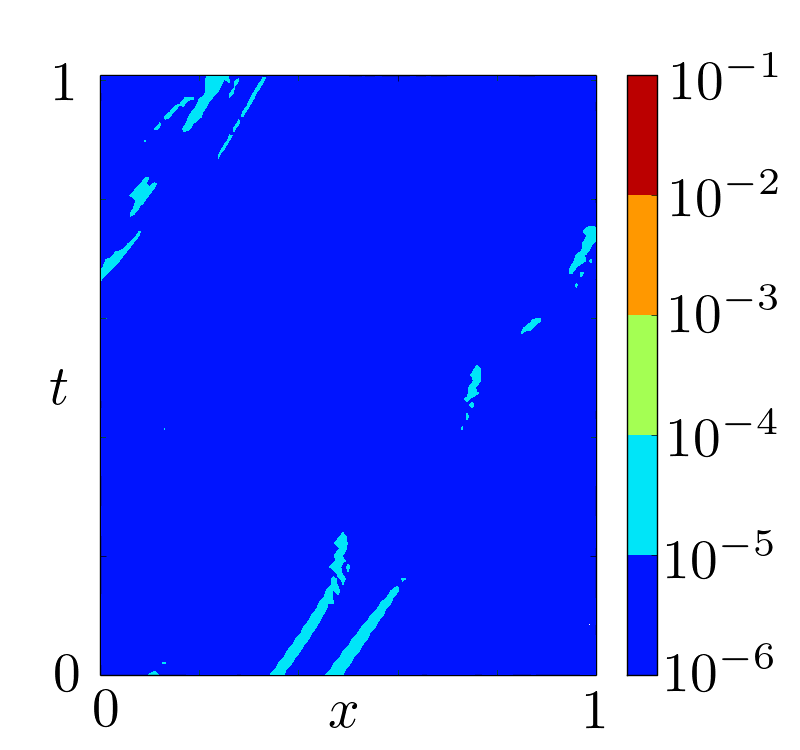}
        \caption{Solution 3}
        \label{fig: grad BL 3}
    \end{subfigure}
    \caption{The errors of estimated gradients for the three solutions.}
    \label{fig: adap basis grad err BL}
\end{center}\end{figure}

Consider a time-space-dependent control $c=c(t,x)$ in \eqref{eqn: Buckley-Leverett} and 
\eqref{eqn: Buckley-Leverett twin}. The gradient of $\xi$, \eqref{eqn: objective ad hoc}, is
estimated using the trained twin model. The estimated gradients are compared with the true adjoint
gradients of the gray-box model, and the errors are shown in Figure \ref{fig: adap basis grad err BL}.\\

The adaptive basis construction improves the accuracy of the gradient estimation.
Table \ref{table: BL grad error} shows the integrated gradient error 
\footnote{The gradient error is integrated using the same quadrature rule as in 
\eqref{eqn: solution mismatch}.} by using either the
ad hoc bases in Figure \ref{fig: sigmoid basis ad hoc} or by using the bases constructed adaptively.\\
\begin{center}
    \begin{tabular}{cccc}
       \hline
         & Solution 1 & Solution 2 & Solution 3\\
       \hline
       Ad hoc basis & $2.5\times 10^{-3}$& $6.6\times 10^{-4}$ & $7.3\times 10^{-5}$ \\
       \hline
       Adaptive basis & $4.2\times 10^{-6}$& $1.5\times 10^{-6}$ & $8.9\times 10^{-7}$ \\
       \hline
    \end{tabular}
    \captionof{table}{The integrated errors of the estimated gradients for the three solutions.}
    \label{table: BL grad error}
\end{center}

\subsection{Navier-Stokes Flow}
\label{sec: chap2 num example NS}
\begin{figure}\begin{center}
    \includegraphics[width=8cm]{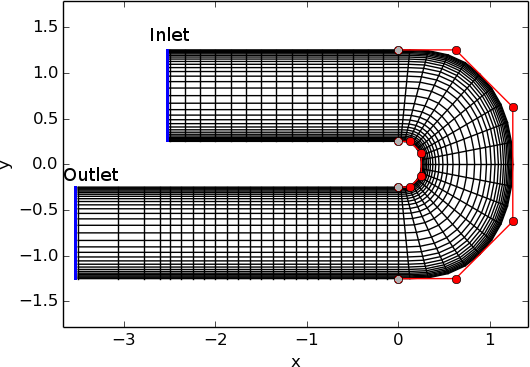}
    \caption{The return bend geometry and the mesh for the simulation.}
    \label{fig: NS mesh}
\end{center}\end{figure}

Consider
a compressible internal flow in a 2-D return bend channel 
driven by the pressure difference between the inlet and the outlet.
The return bend is bounded by 
no-slip walls. The inlet static pressure and the outlet pressure are fixed.
The geometry of the return bend is given in Figure \ref{fig: NS mesh}.
The inner and outer boundaries of the bending section are each generated by 6 control points
using quadratic B-spline.\\

The flow is governed by Navier-Stokes equations.
Let $\rho$, $u$, $v$, $E$, and $p$ denote the density, Cartesian velocity components, 
total energy, and pressure.
The steady-state Navier-Stokes equation is
\begin{equation}
    \frac{\partial}{\partial x} 
    \begin{pmatrix}
        \rho u\\
        \rho u^2 + p - \sigma_{xx}\\
        \rho uv - \sigma_{xy}\\
        u(E\rho+p) - \sigma_{xx} u - \sigma_{xy} v
    \end{pmatrix}
    + \frac{\partial}{\partial y}
    \begin{pmatrix}
        \rho v\\
        \rho uv-\sigma_{xy}\\
        \rho v^2+p-\sigma_{yy}\\
        v(E\rho+p) - \sigma_{xy} u -\sigma_{yy}v
    \end{pmatrix} 
    = \boldsymbol{0}\,,
    \label{NSeqn}
\end{equation}
where
\begin{equation}\begin{split}
    \sigma_{xx} &= \mu \left(2 \frac{\partial u}{\partial x} - \frac{2}{3} \left(\frac{\partial u}{\partial x} 
    + \frac{\partial v}{\partial y}\right)\right)\\
    \sigma_{yy} &= \mu \left(2 \frac{\partial v}{\partial y} - \frac{2}{3} \left(\frac{\partial u}{\partial x} 
    + \frac{\partial v}{\partial y}\right)\right)\\
    \sigma_{xy}&=\mu\left(\frac{\partial u}{\partial y} + \frac{\partial v}{\partial x}\right)
\end{split}\,.\end{equation}
The Navier-Stokes equation requires an additional equation, the state equation, for closure.
The state equation has the form
\begin{equation}
    p = p(U, \rho)\,,
    \label{state equation}
\end{equation}
where $U$ denotes the internal energy per unit volume,
\begin{equation}
    U = \rho\left(E-\frac{1}{2}(u^2+v^2)\right)\,.
\end{equation}

Many models have been developed for the state equation, such as the ideal gas equation, the
van der Waals equation, and the Redlich-Kwong equation \cite{state eqns}.
In the sequel, the true state equation in the gray-box simulator is assumed unknown and
will be inferred from the gray-box solution. 
Let $\rho_\infty$ be the steady state density,
$\boldsymbol{u}_\infty = (u_\infty, v_\infty)$ be the steady state Cartesian velocity,
and $E_\infty$ be the steady state energy density.
The steady state mass flux is
\begin{equation}
    \xi = - \int_{\textrm{outlet}} \rho_\infty u_\infty \big|_{\textrm{outlet}} \; dy=
    \int_{\textrm{inlet}} \rho_\infty u_\infty\big|_{\textrm{inlet}} \; dy
    \label{eqn: mass flux}
\end{equation}
The goal is to estimate the gradient of $\xi$
to the red control points' coordinates.
\\

Two state equations are tested: the ideal gas equation and the Redlich-Kwong equation, 
given by
\begin{equation}\begin{split}
    p_{ig} &= (\gamma-1) U\\
    p_{rk} &= \frac{(\gamma-1)U}{1-b_{rk}\rho} - 
    \frac{a_{rk}\rho^{5/2}}{((\gamma-1)U)^{1/2}(1+b_{rk}\rho)}
\end{split}\label{NS state equations}
\end{equation}
where $a_{rk}=10^7$ and $b_{rk}=0.1$.\\

\begin{figure}\begin{center}
    \begin{subfigure}[t]{.49\textwidth}
        \centering
        \includegraphics[height=8cm]{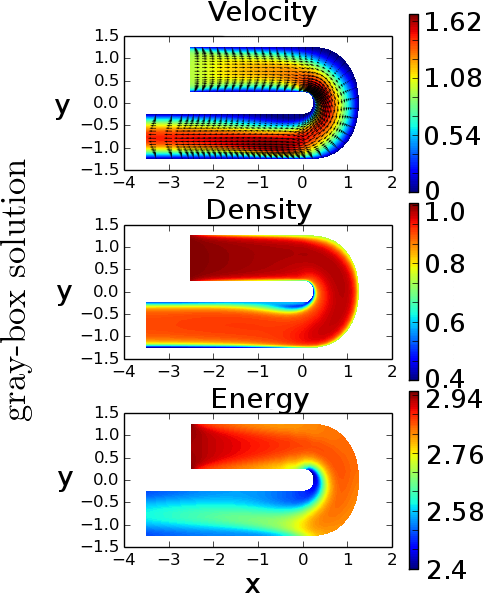}
        \label{fig: graysol Ubend}
    \end{subfigure}
    \begin{subfigure}[t]{.49\textwidth}
        \centering
        \includegraphics[height=8cm]{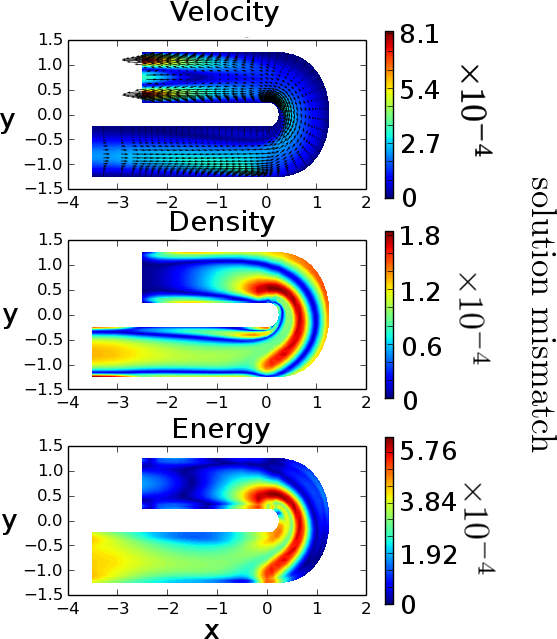}
        \label{fig: errsol Ubend}
    \end{subfigure}
    \caption{Left column: an example gray-box solution for a given geometry. Right column:
             the solution mismatch after training a twin model.}
    \label{fig: grayErrSol Ubend}
\end{center}\end{figure}

\begin{figure}\begin{center}
    \begin{subfigure}[t]{.99\textwidth}
        \centering
        \includegraphics[height=4.5cm]{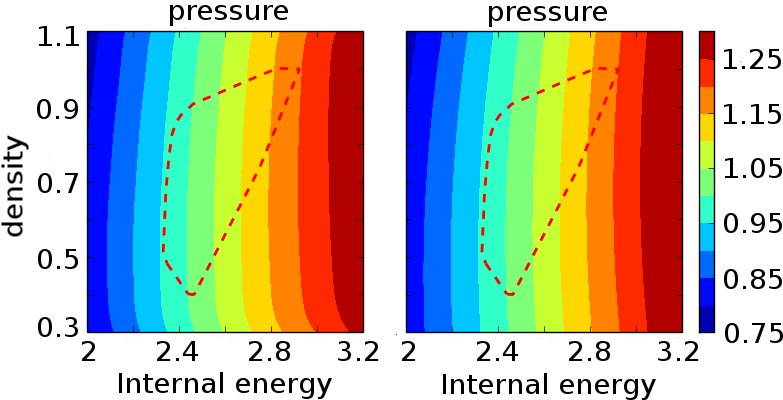}
        \label{fig: graysol Ubend}
    \end{subfigure}\\
    \begin{subfigure}[t]{.99\textwidth}
        \centering
        \includegraphics[height=4.5cm]{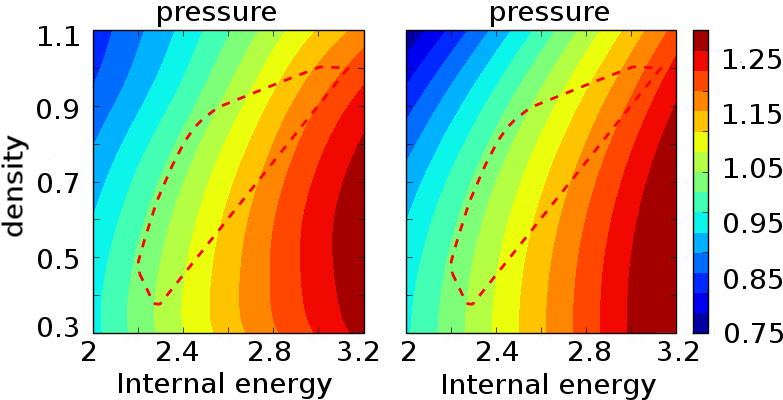}
        \label{fig: errsol Ubend}
    \end{subfigure}
    \caption{The gray-box state equation (right column) and the trained state equation 
             (left column). The gray-box model uses either the ideal gas equation
             (first row) or the Reclich-Kwong equation (second row). The convex hull 
             of the gray-box solution is shown by the dashed red line.}
    \label{fig: state err Ubend}
\end{center}\end{figure}

The solution mismatch, \eqref{eqn: solution mismatch}, is given by
\begin{equation*}\begin{split}
    \mathcal{M} = &w_\rho \int_\Omega \left|\tilde{\rho}_{\infty} - \rho_{\infty}\right|^2 d\boldsymbol{x}
                + w_u
                \int_\Omega \left|\tilde{u}_{\infty}- u_{\infty}\right|^2 d\boldsymbol{x} \\
                + &w_v
                \int_\Omega \left| \tilde{v}_{\infty}- v_{\infty}\right|^2 d\boldsymbol{x}
                + w_E
                \int_\Omega \left|\tilde{E}_{\infty} - E_\infty\right|^2 d\boldsymbol{x}\,,
    \label{NS mismatch}
\end{split}\end{equation*}
where $w_\rho$, $w_u$, $w_v$, and $w_E$ are non-dimensionalization constants.
Figure \ref{fig: grayErrSol Ubend} shows the gray-box solution and the solution mismatch
after training the twin model \footnote{Both the solution and the mismatch are normalized.}.
Figure \ref{fig: state err Ubend} compares the true state equation and the corresponding trained 
state equation, where the convex hull of $({U}_\infty, \rho_\infty)$, the internal energy
and the density of the gray-box solution,
is shown by the dashed red line.
Because the state equation is expected to be inferrable only inside the domain of the gray-box 
solution, large deviation is expected outside the convex hull.\\

\begin{figure}\begin{center}
    \begin{subfigure}[t]{.45\textwidth}
        \centering
        \includegraphics[height=6.2cm]{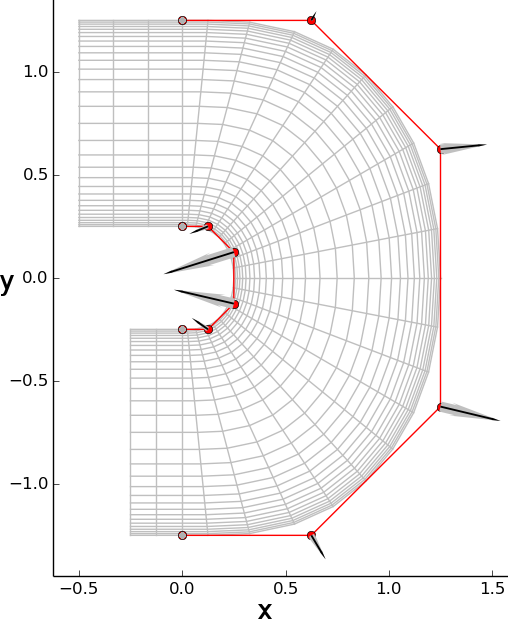}
        \caption{
        The gradient of $\xi$ to the control points for the 
        Redlich-Kwong gas. 
        The wide gray arrow is the gradient evaluated by the gray-box model, while
        the thin black arrow is the gradient evaluated by the twin model using finite difference.
        }
        \label{fig: geo grad}
    \end{subfigure}
    \hspace{.5cm}
    \begin{subfigure}[t]{.45\textwidth}
        \centering
        \includegraphics[height=6.2cm]{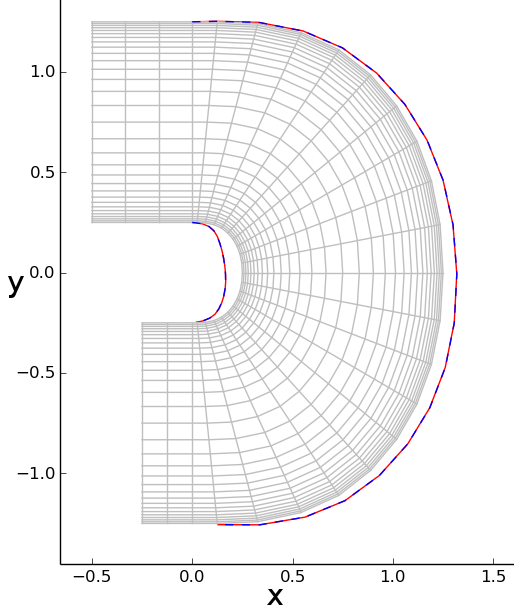}
        \caption{
        The boundary perturbed according to the gradient. 
        The blue dashed line is computed by finite difference of the gray-box model, while the
        red dashed line is computed by the twin model's gradient.
        }
        \label{fig: geo grad perturb}
    \end{subfigure}
    \caption{A comparison of the estimated gradient and the true gradient.}
    \label{fig: geo grad all}
\end{center}\end{figure}

The trained twin model enables the adjoint gradient estimation.
Figure \ref{fig: geo grad all} shows the
estimated gradient of $\xi$ with respect to the control points coordinates. It also
compares the estimated gradient with the true gradient. The two gradients are indistinguishable,
and the error is given in Table \ref{tab: idea gas gradient}.

\begin{center}
\begin{tabular}{ccC{12mm}C{12mm}C{12mm}C{12mm}C{12mm}C{12mm}C{12mm}C{12mm}}
\hline
{Gas}&\multicolumn{4}{c}{Interior control points}&\multicolumn{4}{c}{Exterior control points}\\
\cline{1-9}
Ideal &0.13 & 0.04 & 0.05 & 0.32 &0.16 & 0.15 & 0.07 & 0.02\\
\cline{1-9}
Redlich-Kwong& 0.32 & 0.03 & 0.07 & 0.50 & 0.40 & 0.12 & 0.06 & 0.05\\
\hline
\end{tabular}
\captionof{table}{The error of the gradient estimation, in percentage.}
\label{tab: idea gas gradient}
\end{center}

\subsection{Polymer Injection in Petroleum Reservoir}
\label{sec: chap 2 reservoir}
\begin{figure}
    \begin{center}
        \includegraphics[width=6cm]{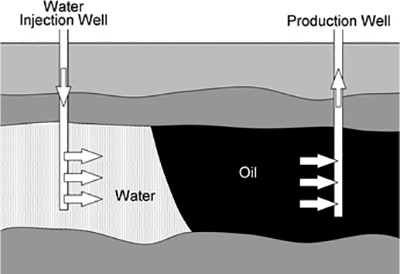}
        \caption{Water flooding in petroleum reservoir engineering (courtesy
        from PetroWiki).}
        \label{fig: polymer sketch}
    \end{center}
\end{figure}

Water flooding is a technique to enhance the secondary recovery in petroleum reservoirs, as illustrated
in Figure \ref{fig: polymer sketch}. 
Injecting pure water can be cost-inefficient due to low water viscosity
and high water cut. Therefore, water-solvent polymer can be utilized to increase the water-phase
viscosity and to reduce the residual oil.\\

Consider a reservoir governed by the two-phase porous media flow equations
\begin{equation}\begin{split}
    \frac{\partial }{\partial t} \left(\rho_\alpha \phi S_\alpha \right) + \nabla \cdot
    \left( \rho_\alpha \vec{v}_{\alpha} \right) &= 0\,, \quad \alpha \in \{w,o\}\\
    \frac{\partial}{\partial t}\left( \rho_w \phi S_w c \right) + \nabla \cdot
    \left( c \rho \vec{v}_{wp}\right) &= 0        
    \end{split}\,,
    \label{eqn: two phase polymer}
\end{equation}
for $x\in \Omega$ and $t\in [0,T]$,
where the phase velocities are given by the Darcy's law
\begin{equation}\begin{split}
    \vec{v}_\alpha &= - {M_\alpha} k_{r\alpha} \boldsymbol{K} \cdot (\nabla p - \rho_w g \nabla z), \, \quad \alpha \in \{w,o\}\\
    \vec{v}_{wp} &= -{M_{wp}} k_{rw} \boldsymbol{K} \cdot (\nabla p - \rho_{w} g \nabla z)
\end{split}\,.
\label{eqn: darcy law}
\end{equation}
$w, o$ indicate the water and oil phases.
$\rho$ is the phase density. $\phi$ is the porosity. $S$ is the phase saturation where
$S_w+S_o=1$.
$c$ is the polymer concentration in the water phase. $v_{w}$, $v_{o}$, 
$v_{wp}$ are the componentwise velocities of water, oil, and polymer. 
$\boldsymbol{K}$
is the permeability tensor. $k_{r}$ is the relative permeability. $p$ is the pressure. $z$ is the depth.
$g$ is the gravity constant. The mobility factors, $M_o, M_w, M_{wp}$, 
model the modification of the componentwise mobility due to the presence of polymer.
In the sequel, the models for the mobility factors are unknown. The only 
knowledge about the mobility factors is that they depend on $S_w, p$, and $c$.\\

\emph{PSim}, the simulator aforementioned in Section \ref{sec: motivation}, is
used as the gray-box simulator, which uses the IMPES time marching, i.e. implicit in pressure
and explicit in saturation, as well as the upwind scheme.
Its solution, $S_w$, $c$, and $p$ can be used to train the 
twin model. The twin model uses fully implicit time marching and the upwind scheme.
The solution mismatch is defined by
\begin{equation}
    \mathcal{M} = w_{S_w}\int_0^T\int_\Omega |S_w-\tilde{S}_w|^2 d\boldsymbol{x} dt
                + w_{c}\int_0^T\int_\Omega |c-\tilde{c}|^2 d\boldsymbol{x} dt 
                + w_{p}\int_{0}^T\int_\Omega |p-\tilde{p}|^2 d\boldsymbol{x} dt\,,
    \label{eqn: polymer sol mismatch}
\end{equation}
where $w_{S_w}$, $w_c$, and $w_p$ are non-dimensionalization constants.\\

\begin{figure}
    \begin{center}
        \includegraphics[height=4.cm]{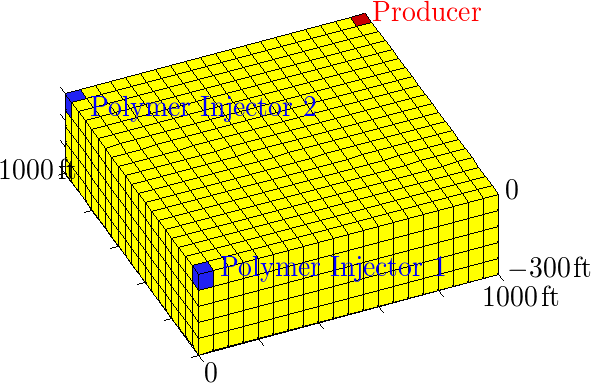}
        \caption{The geometry of the petroleum reservoir.}
        \label{fig: reservoir 3D mesh}
    \end{center}
\end{figure}
\begin{figure}\begin{center}
    \begin{subfigure}[t]{.99\textwidth}
        \centering
        \includegraphics[height=4.cm]{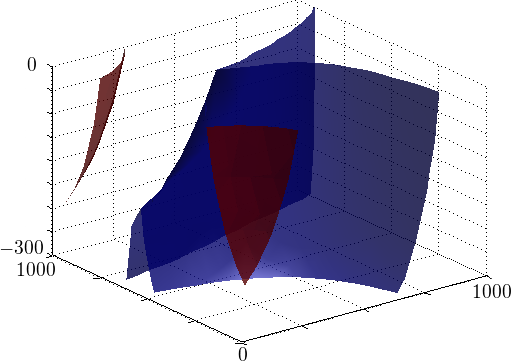}
        \label{fig: reservoir 3D untrained}
        \caption{Untrained twin model.}
    \end{subfigure}\\
    \begin{subfigure}[t]{.99\textwidth}
        \centering
        \includegraphics[height=4.cm]{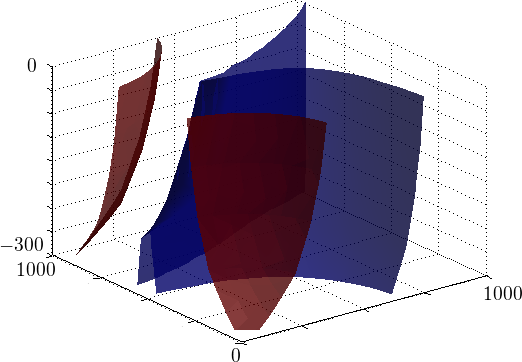}
        \label{fig: reservoir 3D untrained}
        \caption{PSim.}
    \end{subfigure}\\
    \begin{subfigure}[t]{.99\textwidth}
        \centering
        \includegraphics[height=4.cm]{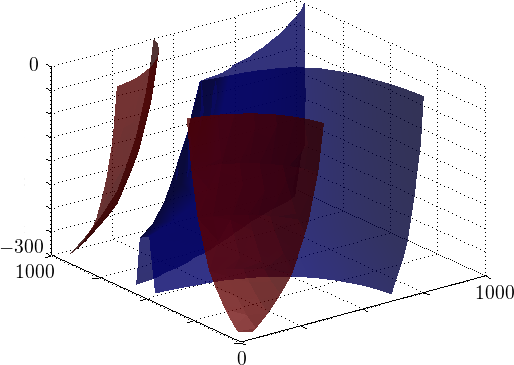}
        \label{fig: reservoir 3D untrained}
        \caption{Trained twin model.}
    \end{subfigure}
    \caption{The isosurfaces of $S_w=0.25$ and $S_w=0.7$ at $t=30$ days. }
    \label{fig: reservoir 3D solutions}
\end{center}\end{figure}

Consider a reservoir setup shown in Figure \ref{fig: reservoir 3D mesh}, which
is a 3D block with two injectors and one producer. The permeability
is 100 milli Darcy, and the porosity is 0.3. A constant injection rate of $10^6 
\texttt{ft}^3/\texttt{day}$ is used at both the injectors.
The reservoir is simulated for $t\in [0,50] \texttt{day}$.
The solution of $S_w$ is illustrated in Figure \ref{fig: reservoir 3D solutions} for the untrained 
twin model, the gray-box model, and the trained twin model, respectively. 
After the training, the twin-model solution matches the gray-box solution closely.\\

\begin{figure}
    \begin{center}
        \includegraphics[width=6.5cm]{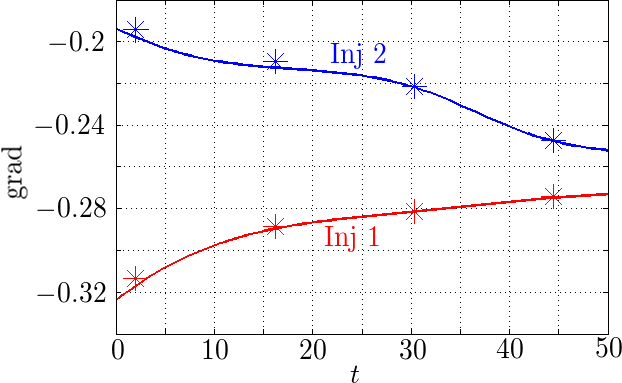}
        \caption{The gradient of $\xi$ with respect to rates at the two injectors.
                 The lines indicate the gradients estimated by the twin model, while
                 the stars indicate the true gradient evaluated by finite difference.}
        \label{fig: reservoir 3D gradient}
    \end{center}
\end{figure}

Let the objective function be the residual oil at $T=50 \, \texttt{day}$,
\begin{equation}
    \xi = \int_\Omega \rho_o(T) \phi S_o(T) \,\textrm{d} \boldsymbol{x} \,.
    \label{eqn: chap2 reservoir xi day 50}
\end{equation}
The gradient of $\xi$ with respect to the time-dependent injection rate is computed.
The gradient estimated by the twin model is shown in Figure \ref{fig: reservoir 3D gradient},
where the red and blue lines indicate the gradient for the two injectors. In comparison, 
the star markers show the true gradient at day $2$, $16$, $30$, and $44$, evaluated by finite difference.
Clearly, a rate increase at the injector 1 leads to more residual oil reduction than the injector 2.
This is because the injector 2 is closer to the producer, where a larger rate accelerates the 
water breakthrough that impedes further oil production.
It is observed that the estimated gradient closely matches the true gradient, although
the error slightly increases for smaller $t$, possibly because of the different numerical schemes
used in the twin and gray-box models. The error is given in Table \ref{tab: reservoir 3D grad error}.\\

\begin{center}
    \begin{tabular}{ccccc}
       \hline
         Error & $t=0.04$ & $t=0.32$ & $t=0.6$ & $t=0.88$ \\ \hline
         Inj 1 & 1.7  & 1.0 & 0.6 & 0.2  \\ \hline
         Inj 2 & 2.2 &  1.9 & 0.7 & 0.2  \\ \hline
    \end{tabular}
    \captionof{table}{The error of estimated gradient at day $2$, $16$, $30$, and $44$, in percentage.}
    \label{tab: reservoir 3D grad error}
\end{center}

\section{Conclusions}
\label{sec: chap 2 summary}
This article develops a method for gradient estimation by using the space-time solution of
gray-box conservation law simulations. In particular, an adjoint-enabled twin model 
is trained to minimize the solution mismatch metric.
The inferrability of the twin model is studied theoretically for a simple PDE with only one equation
and one dimensional space. 
To enable the training computationally, a sigmoid parameterization
is presented. However, an ad hoc choice for the bases
does not fully exploit the information contained in the
gray-box solution. To address this issue,
an adaptive basis construction procedure is presented. The adaptive 
procedure  builds upon three key elements: the approximated basis significance,
the basis neighborhood, and the cross validation. The algorithm for training the twin model
is summarized. To alleviate the training cost, a pre-train step is suggested
that minimizes the integrated truncation error instead of the solution mismatch.\\

The proposed twin model algorithm has a wide applicability, which is demonstrated on 
a variety of numerical examples. The first example is the Buckley-Leverett equation,
whose flux function is inferred. The trained twin model accurately estimates 
the gradient of an objective to the source term. 
The second example is the steady-state Navier-Stokes equation in a return bend,
whose state equation is inferred. The inferred state equation allows 
estimating the gradient of mass flux to the control surface geometry.
The third example is the petroleum reservoir with polymer injection, where
the mobility factors are inferred. The gradient of the residual oil to
the injection rate is estimated.
With the aid of the estimated gradient, the objective can be optimized more efficiently,
which will be discussed in a companion paper.\\

There are several potential thrusts of further research: 
A useful extension is to investigate the
inferrability of twin models for various conservation laws. In particular,
Theorem \ref{theorem: 1} may be extended for problems with a
system of equations and higher spatial dimension. Another interesting extension
is to study the applicability of the pre-train step, especially to
obtain a necessary and sufficient condition for bounding $\mathcal{M}$
with $\mathcal{T}$. Finally, it is interesting 
to generalize the formulation \eqref{eqn: govern PDE} to
incorporate unknown source terms and boundary conditions.\\

\section{Appendix}

\subsection{Theorem \ref{theorem: 1}}
\label{proof 1}
Proof: \\
We prove false the contradiction of the theorem, which reads:\\
\emph{For any $\delta>0$ and $T>0$, there exist $\epsilon > 0$, and 
$F, \tilde{F}$ satisfying the conditions stated
in theorem \ref{theorem: 1}, such that $\|\tilde{u} - u\|_{\infty} < \delta$ and 
$\left\| \frac{d\tilde{F}}{du} - \frac{dF}{du} \right\|_\infty > \epsilon$ on $B_u$.}\\

We show the following exception to the contradiction in order to prove it false.\\
\emph{For any $\epsilon >0$ and any $F,\tilde{F}$ satisfying $\left\|\frac{d\tilde{F}}{du} - \frac{dF}{du}\right\|_{\infty} > \epsilon$
on $B_u$, we can find $\delta >0$ and $T>0$ such that $\|\tilde{u} - u\|_\infty > \delta$.
}\\

The idea is to construct such an exception by the method of lines \cite{method of lines}.
Firstly, assume there is no shock wave for \eqref{eqn: easy 1} and \eqref{eqn: easy 2} for $t\in[0,T]$. 
Choose a segment in space, $[x_0-\Delta, x_0]$ with $0< \Delta< \frac{\epsilon}{L_{F} L_{u}}$, that satisfies 
\begin{itemize}
    \item $u_0(x) \in B_u$ for any $x\in [x_0-\Delta, x_0]$;
    \item $\left| \frac{d\tilde{F}}{du}\big(u_0(x_0)\big) - \frac{dF}{du} \big(u_0(x_0)\big) \right| > \epsilon$;
    \item $x_0 - \Delta + \frac{dF}{du}\big( u_0(x_0 - \Delta) \big) T = x_0 + \frac{d\tilde{F}}{du}\big( u_0 \big)T \equiv x^*$.
\end{itemize}
Without loss of generality, we assume $\frac{dF}{du} >0$ and $\frac{d\tilde{F}}{du}>0$ for $\big\{u \big|u=u_0(x)\,,\; x\in[x_0-\Delta, x_0]\big\}$.
Using the method of lines, we have
$$u\left(T\,,\; x_0-\Delta +\frac{dF}{du}\big( u_0(x_0-\Delta) \big)T  \right) = u_0(x_0 -\Delta)\,,$$ and
$$\tilde{u}\left(T\,,\; x_0+\frac{d\tilde{F}}{du}\big( u_0(x_0) \big)T  \right) = u_0(x_0 )\,.$$
Therefore
$$
    \left| \tilde{u}(x^*, T) - u(x^*, T) \right| = \left| u_0(x_0) - u_0(x_0-\Delta) \right| \ge \gamma \Delta \equiv \delta\,,
$$
by using the definition of $B_u$.\\

Set 
$T= \frac{\Delta}{\left| \frac{d\tilde{F}}{du}\big( u_0(x_0-\Delta)\big) - \frac{dF}{du}\big( u_0(x_0) \big) \right|}$, we have
\begin{equation*}\begin{split}
     & \left| \frac{d\tilde{F}}{du}\big( u_0(x_0-\Delta)\big) - \frac{dF}{du}\big( u_0(x_0) \big) \right| \\
    = &\left| \frac{d{F}}{du}\big( u_0(x_0)\big) - \frac{d\tilde{F}}{du}\big( u_0(x_0) \big) + 
           \frac{dF}{du} \big( u_0(x_0-\Delta) \big) - \frac{dF}{du}\big( u_0(x_0)\big) \right|\\
    = &\left| \frac{d{F}}{du}\big( u_0(x_0)\big) - \frac{d\tilde{F}}{du}\big( u_0(x_0) \big) + 
       \overline{\frac{d^2 F}{du^2}} \big( u_0(x_0-\Delta) - u_0(x_0) \big) \right|\\
    \ge& \left| \frac{dF}{du}\big( u_0(x_0) \big) - \frac{d\tilde{F}}{du}\big( u_0(x_0) \big) \right| - L_uL_F \Delta \\
    \ge& \epsilon - L_uL_F\Delta\equiv \epsilon_{F} > 0
\end{split}\end{equation*}
by using the mean value theorem. Therefore $T \le \frac{\Delta}{\epsilon_F} < \infty$.
So we find a $\delta = \gamma \Delta$ and a $T < \infty$ that provides an exception
to the contradiction of the theorem. 

Secondly, if there is shock wave within $[0,T]$ for either \eqref{eqn: easy 1} or \eqref{eqn: easy 2}, we let
$T^*$ be the time of the shock occurrence. Without loss of generality, assume the shock occurs for \eqref{eqn: easy 1} first.
The shock implies the intersection of two characteristic lines. Choose a $\Delta>0$
such that $\big| \frac{dF}{du}\big(u_0(x)\big) - \frac{dF}{du} \big( u_0(x-\Delta) \big) \big|T^* = \Delta$.
Using the mean value theorem, we have
$$
    T^* = \frac{\Delta}{ \overline{\frac{d^2 F}{du^2}} \big( u_0(x) - u_0(x-\Delta) \big) } \ge \frac{1}{L_u L_F}
$$
Thus, if we choose
$$
    T = \min \left\{ 
        \frac{1}{L_uL_K},\;\frac{\Delta}{\epsilon_\Delta}
    \right\}\,,
$$
no shock occurs in $t\in[0,T]$. Since the theorem is already proven for the no-shock scenario, the proof completes.\hfill $\surd$

\subsection{Theorem \ref{theorem: 2}}
\label{proof 2}
\begin{figure}\begin{center}
    \includegraphics[height=5.5cm]{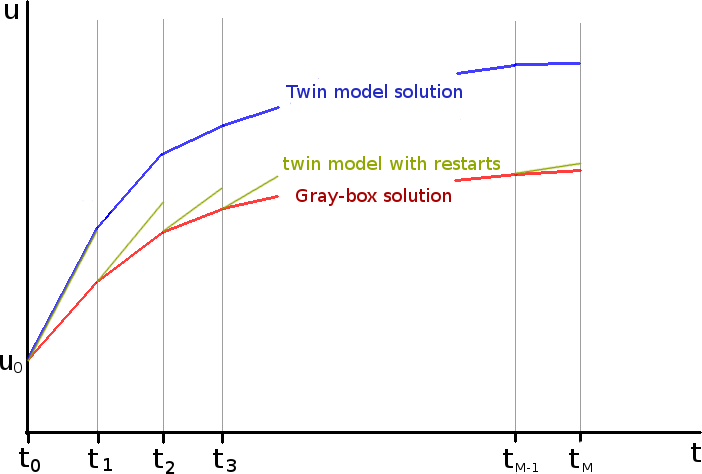}
    \caption{The state-space trajectories of the gray-box model and the twin model. 
             $\mathcal{M}$ measures the difference of the twin model trajectory (blue) with
             the gray-box trajectory (red). $\mathcal{T}$ measures the difference of the
             twin model trajectory with restarts (green) and the gray-box trajectory (red).}
    \label{fig:sketch}
\end{center}\end{figure}

Proof: \\
Let the one-step time marching of the gray-box simulator be
\begin{equation*}
    \mathcal{H}:\, \mathbb{R}^n\mapsto\mathbb{R}^n,\, 
    \boldsymbol{u}_{i\cdot}
    \rightarrow 
    \boldsymbol{u}_{i+1\cdot} = \mathcal{H}_i\boldsymbol{u}_{i\cdot} \,,
    \quad i=1,\cdots, M-1\,,
\end{equation*}
The integrated truncation error can be written as
\begin{equation*}\begin{split}
    \mathcal{T}(\tilde{F}) &= 
        \sum_{i=1}^M \sum_{j=1}^N w_{j} \left(
        \boldsymbol{u}_{i+1\, j} - (\mathcal{G}\boldsymbol{u}_{i\, \cdot})_j
        \right)^2\\
    &= \sum_{i=1}^M (u_{i+1\,\cdot} - \mathcal{G} u_{i\,\cdot})^T W
                    (u_{i+1\,\cdot} - \mathcal{G} u_{i\,\cdot}) \\
    &= \sum_{i=1}^M \left\|u_{i+1\,\cdot} - \mathcal{G} u_{i\,\cdot}\right\|^2_{W}\\
    &= \sum_{i=1}^M \left\| \mathcal{H} u_{i\, \cdot} - \mathcal{G} u_{i\,\cdot} \right\|^2_{W}\\
    &= \sum_{i=1}^M \left\| \left(\mathcal{H}^i - \mathcal{G} \mathcal{H}^{i-1}\right)
       u_{0\,\cdot} \right\|^2_{W}\,.
\end{split}\end{equation*}
Similarly, the solution mismatch can be written as
\begin{equation*}
    \mathcal{M}(\tilde{F}) = \sum_{i=1}^M \left\| 
    \left(\mathcal{H}^i - \mathcal{G}^i\right)
    u_{0\,\cdot} \right\|^2_{W}
\end{equation*}
Fig \ref{fig:sketch} gives an explanation of $\mathcal{M}$ and $\mathcal{T}$ by viewing 
the simulators as discrete-time dynamical systems.\\

Using the equality
\begin{equation*}
    \mathcal{G}^i-\mathcal{H}^i = (\mathcal{G}^i-\mathcal{G}^{i-1}\mathcal{H}) + (\mathcal{G}^{i-1}
    \mathcal{H} - \mathcal{G}^{i-2}\mathcal{H}^2) + \cdots + (\mathcal{G}\mathcal{H}^{i-1}-\mathcal{H}^i)\,,\quad
    i\in \mathbb{N}\,,
\end{equation*}
and triangular inequality, we have 
\footnotesize{
\begin{equation*}
    \mathcal{M} \le 
    \left\{\begin{split}
        \|(\mathcal{G}^{M-1} \mathcal{G} - \mathcal{G}^{M-1}\mathcal{H})u_{0\cdot}\|^2_{W} &+ \|(\mathcal{G}^{M-2}\mathcal{G}\mathcal{H} - \mathcal{G}^{M-2}\mathcal{H}^2)u_{0\cdot}\|^2_{W}&+\cdots
        &+ \|(\mathcal{G}\mathcal{H}^{M-1} - \mathcal{H}^{M})u_{0\cdot}\|^2_{W}\\
        &+\|(\mathcal{G}^{M-2} \mathcal{G} - \mathcal{G}^{M-2}\mathcal{H})u_{0\cdot}\|^2_{W} &+ \cdots 
        &+ \|(\mathcal{G}\mathcal{H}^{M-2} - \mathcal{H}^{M-1})u_{0\cdot}\|^2_{W}\\
        &\ddots&& \vdots\\
        &&& + \|(\mathcal{G} - \mathcal{H})u_{0\cdot}\|^2_{W}
    \end{split}
    \right\}\,.
    \label{expansion global error}
\end{equation*}
}
\normalsize
Therefore,
\begin{equation*}
    \mathcal{M} - \mathcal{T} \le \qquad\hspace{11cm}
\end{equation*}
\footnotesize
\begin{equation*}
    \left\{\begin{split}
        \|(\mathcal{G}^{M-1} \mathcal{G} - \mathcal{G}^{M-1}\mathcal{H})u_{0\cdot}\|^2_{W} &+ \|(\mathcal{G}^{M-2}\mathcal{G}\mathcal{H} - \mathcal{G}^{M-2}\mathcal{H}^2)u_{0\cdot}\|^2_{W}&+\cdots
        &+ \|(\mathcal{G}\mathcal{G}\mathcal{H}^{M-2} - \mathcal{G} \mathcal{H}^{M-1})u_{0\cdot}\|^2_{W}\\
        &+\|(\mathcal{G}^{M-2} \mathcal{G} - \mathcal{G}^{M-2}\mathcal{H})u_{0\cdot}\|^2_{W} &+ \cdots 
        &+ \|(\mathcal{G}\mathcal{G}\mathcal{H}^{M-3} - \mathcal{G} \mathcal{H}^{M-2})u_{0\cdot}\|^2_{W}\\
        &\ddots&& \vdots\\
        &&& + \|(\mathcal{G} \mathcal{G} - \mathcal{G} \mathcal{H})u_{0\cdot}\|^2_{W}
    \end{split}
    \right\}\,.
    \label{error diff}
\end{equation*}
\normalsize
Under the assumption
\begin{equation*}
    \left\| \mathcal{G}a - \mathcal{G}b \right\|^2_{W} \le \beta \|a-b\|^2_{W}\,,
\end{equation*}
and its implication
\begin{equation*}
    \left\|\mathcal{G}^i a - \mathcal{G}^i b\right\|^2_W \le \beta^i \|a-b\|_W^2\,, \quad i\in \mathbb{N}\,,
\end{equation*}
we have
\begin{equation*}
    \mathcal{M} - \mathcal{T} \le 
    \left\{\begin{split}
        \beta^{M-1}\|(\mathcal{G} - \mathcal{H})u_{0\cdot}\|^2_W &+ \beta^{M-2}\|(\mathcal{G}\mathcal{H} - \mathcal{H}^2)u_{0\cdot}\|^2_W&+\cdots
        &+ \beta\|(\mathcal{G}\mathcal{H}^{M-2} - \mathcal{H}^{M-1})u_{0\cdot}\|^2_W\\
        &+\beta^{M-2}\|( \mathcal{G} - \mathcal{H})u_{0\cdot}\|^2_{M-1} &+ \cdots 
        &+ \beta\|(\mathcal{G}\mathcal{H}^{n-3} - \mathcal{H}^{n-2})u_{0\cdot}\|^2_W\\
        &\ddots&& \vdots\\
        &&& + \beta\|(\mathcal{G} - \mathcal{H})u_{0\cdot}\|^2_W
    \end{split}
    \right\}\,.
    \label{error diff}
\end{equation*}
Reorder the summation, we get
\begin{equation*}
    \mathcal{M} - \mathcal{T} \le 
    \left\{\begin{split}
        \beta^{M-1}\|(\mathcal{G} - \mathcal{H})u_{0\cdot}\|^2_W &+ \beta^{M-2}\|(\mathcal{G} - \mathcal{H})u_{0\cdot}\|^2_W&+\cdots
        &+ \beta\|(\mathcal{G} - \mathcal{H})u_{0\cdot}\|^2_W\\
        &+\beta^{M-2}\|( \mathcal{G}\mathcal{H} - \mathcal{H}^2)u_{0\cdot}\|^2_W &+ \cdots 
        &+ \beta\|(\mathcal{G}\mathcal{H} - \mathcal{H}^{2})u_{0\cdot}\|^2_W\\
        &\ddots&& \vdots\\
        &&& + \beta\|(\mathcal{G}\mathcal{H}^{M-2} - \mathcal{H}^{M-1})u_{0\cdot}\|^2_W
    \end{split}
    \right\}\,.
    \label{error diff 2}
\end{equation*}
Therefore,
\begin{equation*}
    \mathcal{M} - \mathcal{T} \le \left(\beta^{M-1} + \beta^{M-2} + \cdots + \beta \right)
    \mathcal{T}
\end{equation*}
If $\beta$ is strictly less than $1$, then
\begin{equation*}
    \mathcal{M} \le \frac{1}{1-\beta} \mathcal{T}\,,
\end{equation*}
thus completes the proof.\hfill $\surd$

\clearpage
\bibliographystyle{ctr}
\bibliography{ctr_summer_lesopt}

\end{document}